\documentclass[12pt]{article}
\usepackage{latexsym, amsmath, amsfonts, amssymb}
\usepackage{hyperref}
\newcommand{\nc}{\newcommand}
\nc{\ga}{\gamma} \nc{\di}{\displaystyle}
\nc{\ek}{\protect\\[1ex]}
\nc{\C}{{\mathbb C}} \nc{\N}{{\mathbb N}} \nc{\R}{{\mathbb R}}
\nc{\Z}{{\mathbb Z}} \nc{\La}{\Lambda} \nc{\la}{\lambda}
\nc{\da}{\delta} \nc{\Da}{\Delta} \nc{\na}{\nabla} \nc{\vp}{\varphi}
\nc{\ka}{\kappa} \nc{\si}{\sigma} \nc{\Si}{\Sigma} \nc{\al}{\alpha}
\nc{\be}{\beta} \nc{\om}{\omega} \nc{\Om}{\Omega} \nc{\pa}{\partial}
\nc{\ti}{\times}
 \nc{\ve}{\varepsilon} \nc{\ra}{\rightarrow} \nc{\Ra}{\Rightarrow}
\nc{\ran}{\rangle} \nc{\lan}{\langle}
 \nc{\eq}[1]{\mbox{\rm{(\ref{E#1})}}}
 \nc{\EL}{\cal L}
 \nc{\ha}{\frac{1}{2}}
 \nc{\wrar}{\rightharpoonup}
\nc{\hra}{\hookrightarrow} \nc{\supp}{{\rm supp}\,}
\nc{\dense}{\hra^{\hspace{-3mm}d\,\,}} \nc{\curl}{\text{curl}\,}
\renewcommand{\div}{\text{div}\,} \nc{\Id}{\text{Id}}
 \nc{\qed}{\mbox{}\nolinebreak\hfill \rule{2mm}{2mm}}
\numberwithin{equation}{section}

\newtheorem{lemma}{Lemma}[section]
\newtheorem{theorem}[lemma]{Theorem}

\newtheorem{proposition}[lemma]{Proposition}
\newtheorem{remark}[lemma]{Remark}
\newtheorem{definition}{Definition}[section]

\begin{document}
\bibliographystyle{alpha}

\title{Existence of incompressible and immiscible flows in
critical function spaces on bounded domains}
\author{Myong-Hwan Ri$^\dag$, Ping Zhang$^\ddag$}
\date{}
\maketitle

\begin{abstract}
\noindent We study global existence and uniqueness of solutions to
instationary inhomogeneous Navier-Stokes equations on bounded
domains of $\R^n, n\geq 3$, with initial velocity in
$B^0_{q,\infty}(\Om)$, $q\geq n$, and piecewise constant initial
density.
\par
To this end, first, existence for momentum equations with prescribed
density is obtained based on maximal $L^\infty_\ga$-regularity of
the Stokes operator in little Nicolskii space
$b^{s}_{q,\infty}(\Om)$, $s\in\R$, exploited in \cite{RiZh14} and
existence for divergence problem in $b^{-s}_{q,\infty}(\Om)$, $s>0$.
Then, we obtain an existence result for transport equations in the
space of pointwise multipliers for $b^{-s}_{q,\infty}(\Om)$, $s>0$.
Finally, the existence of the inhomogeneous Navier-Stokes equations
is proved via an iterate scheme while the proof of uniqueness is
done via a Lagrangian approach based on the prior results on
momentum equations and transport equation.
\end{abstract}

{\small{\bf 2000 Mathematical Subject Classification:} 35Q30; 35B35;
76D03; 76D07; 76E99 

\let\thefootnote\relax\footnote{\hspace{-0.3cm}
\dag: Institute of Mathematics, State Academy of Sciences, DPR Korea.

\par
\ddag: Academy of Mathematics $\&$ Systems Science and Hua Loo-Keng Center for Mathematical Sciences, The Chinese Academy of
Sciences, Beijing 100190, CHINA, and School of Mathematical Sciences, University of Chinese Academy of Sciences, Beijing 100049, CHINA. 
Email:{zp@amss.ac.cn}
}


\section {Introduction and Main Result}

\par
In this paper, we consider the initial boundary value problem for inhomogeneous
Navier-Stokes equation
\begin{equation}
\label{E1.1}
\begin{array}{rcl}
     \rho_t+\div (\rho u)=0\,\, &\text{in}&(0,T)\ti\Om,\ek
     \rho u_t-\mu\Da u  + \rho (u \cdot \na) u+\na P  = 0 \,\, &\text{in}&(0,T)\ti\Om,\ek
      \div u = 0 \,\, &\text{in} &(0,T)\ti\Om,\ek
     u= 0\, &\text{on}&(0,T)\ti\pa\Om,\ek
     \rho(0,x)=\rho_0,\quad u(0,x)=u_0 &\text{in}&\Om,
\end{array}
\end{equation}
where $0<T\leq \infty$, $\Om\subset \R^n, n\geq 3,$ is a bounded
domain of $C^2$-class, $\rho$ is density of the fluid, $\mu$ is the
dynamic viscosity, and $u$, $P$ are, respectively,
 the velocity and pressure.

The system \eq{1.1} describes the motion of viscous incompressible
flows with variable density. In particular, the motion of a mixture
of different immiscible and incompressible fluids is modeled by
\eq{1.1}.

Kazhikhov \cite{Kaz74} considered \eq{1.1} in whole $\R^n$, $n=2,3$,
to prove existence of a weak solution in energy space and, in
addition, strong solution for small data.

Ladyzhenskaya and Solonnikov considered
 unique solvability
for \eq{1.1} in bounded domain of $\R^n$, $n=2,3$, in \cite{LaSo75},
where existence and uniqueness of solutions are proved when
$\rho_0\in C^1(\bar\Om)$ is positive away from $0$ and the norm of
$u_0$ in $W^{2-2/q,q}(\Om)$, $q>n,$ is small enough. Similar result
has been obtained by Danchin \cite{Da06}, assuming less regularity
on initial data.

If one considers the case where variable density is close to a constant (say $\bar\rho$)
and
writes $\frac{\bar\rho}{\rho}=a+1$, then it is easily checked that $(a, u)$
solves the system:

\begin{equation}
\label{E1.6}
\begin{array}{rcl}
     a_t+(u\cdot \na )a=0\,\, &\text{in}&(0,T)\ti\Om,\ek
     u_t+(u \cdot \na ) u+(1+a)(-\nu\Da u +\na p)  = 0 \,\, &\text{in}&(0,T)\ti\Om,\ek
      \div u = 0 \,\, &\text{in} &(0,T)\ti\Om,\ek
     u= 0\, &\text{on}&(0,T)\ti\pa\Om,\ek
     a(0,x)=a_0,\quad u(0,x)=u_0 &\text{in}&\Om,
\end{array}
\end{equation}
where and in what follows
$$\nu\equiv \frac{\mu}{\bar\rho}\quad\text{and}\quad p\equiv
\frac{P}{\bar\rho}.$$

   Since the system \eq{1.6} is invariant under the scaling
\begin{equation}
\label{E1.5}
\begin{array}{l}
a_\la(t,x)=a(\la^2 t, \la x), u_\la(t,x)=\la u(\la^2 t, \la x),
 \,\, p_\la(x,t)=\la^2 p(\la^2 t, \la x),\ek
{a_0}_\la=a_0(\la x),\,\,
{u_0}_\la=\la u_0(\la x),\quad \la>0,
\end{array}
\end{equation}
to show the existence of solutions to \eq{1.6} in {\it critical
spaces}, i.e., spaces with norms invariant under the scaling
\eq{1.5} is very important.

In \cite{Da03}, Danchin considered unique solvability in some
scaling invariant homogeneous Besov spaces for \eq{1.6} in the whole
space case with $n\geq 2$; more precisely, he showed that if $(a_0,
u_0)\in (\dot{B}^{n/2}_{2,r}(\R^n)\cap L^\infty(\R^n))\ti
\dot{B}^{n/2-1}_{2,1}(\R^n)$ ($r=\infty$ for $n=3$ and $r=1$ for
$n=2$) and
$\|a_0\|_{\dot{B}^{n/2}_{2,\infty}}+\|u_0\|_{\dot{B}^{n/2-1}_{2,1}}$
is small enough,
 then \eq{1.6} has a unique solution $(\rho,u)$ such that
 $$\begin{array}{l}
 a\in BC([0,T),\dot{B}^{n/2}_{2,r}(\R^n))\cap L^\infty(0,T; L^\infty(\R^n)),\ek
 u\in BC([0,T),\dot{B}^{n/2-1}_{2,1}(\R^n))\cap L^1(0,T; \dot{B}^{n/2+1}_{2,1}(\R^n)).
 \end{array}$$
This result was generalized by Abidi \cite{Abi07}
to the case where the space
$(\dot{B}^{n/2}_{2,\infty}(\R^n)\cap L^\infty(\R^n))\ti \dot{B}^{n/2-1}_{2,1}(\R^n)$
for $(a,u)$ in \cite{Da03} is changed to
$\dot{B}^{n/q}_{q,1}(\R^n)\ti \dot{B}^{-1+n/q}_{q,1}(\R^n)$, $1<q<2n$, thus
showing existence for $1<q<2n$ and uniqueness for $1<q\leq n$; the gap in the uniqueness
for $n<q<2n$ was filled by Danchin and Mucha \cite{DaMu12} via Lagrangian approach.
Furthermore, the smallness assumption for initial density variation $a_0$
was relaxed by Abidi, Gui and Zhang in \cite{AGZ12}, \cite{AGZ13}.
In \cite{DaMu09}, Danchin and Mucha proved global well-posedness of \eq{1.1} in half
space $\R^n_+$ under the assumption that $\rho_0$ is close enough to a constant in
$L^\infty(\R^n_+)\cap \dot{W}^1_n(\R^n_+)$ and the norm of
$u_0\in \dot{B}^0_{n,1}(\R^n_+)$ is small enough.

On the other hand, very recently, an attempt to study
\eq{1.6} (equivalently \eq{1.1}) in critical function
spaces, just assuming
initial density merely bounded positively
away from $0$,  has been made, cf. \cite{DaMu13},  \cite{DaZh13}, \cite{HPZ13}.
 Huang, Paicu and Zhang \cite{HPZ13} proved existence of a global solution
 to \eq{1.6} in whole $\R^n$ under a smallness condition
 of $a_0\in L^\infty(\R^n)$ and $u_0\in \dot{B}^{-1+n/q}_{q,r}(\R^n)$,
 $q\in (1,n), r\in (1,\infty)$, and uniqueness of a solution under a slightly
 higher regularity assumption on initial velocity $u_0$;
this result is extended to the half-space setting by Danchin and
Zhang \cite{DaZh13}. The main ideas of \cite{HPZ13} and
\cite{DaZh13} are to employ, for existence, maximal $L^p$-regularity
for Stokes operator in Lebesgue spaces and, for uniqueness, a
Lagrangian approach which was exploited in \cite{DaMu12}. In
\cite{DaMu13}, the case of bounded domains with $C^2$-boundary was
considered under the same assumptions on the initial density and
$u_0\in B^{2-2/r}_{q,r}(\Om)$, $n<q<\infty$, $1<r<\infty$,
 $1/q\neq 2-2/r$.

Here, we recall the famous result by Kato \cite{Kat84} for strong
solvability for classical homogeneous  Navier-Stokes system
($\rho\equiv const$) with initial velocity
\begin{equation}
\label{E1.7}
u_0\in L^n(\R^n),\quad \div u_0=0.
\end{equation}
Natural question whether the initial condition \eq{1.7} will still
guarantee local or global well-posedness for inhomogeneous
Navier-Stokes system \eq{1.1} arises. However, by the above
mentioned previous results for inhomogeneous systems \eq{1.1} or
\eq{1.6}, $u_0$ is allowed to take in $\dot{B}^0_{n,1}(\R^n),
\dot{B}^0_{n,1}(\R^n_+)$,
 or some homo- or inhomogeneous Besov spaces
which, at least, do not include $L^n$-space.
\par\vspace{0.2cm}
Therefore, in this paper, we show the existence of a solution for
\eq{1.1} on bounded domains when initial velocity is in
inhomogeneous Besov spaces $B^0_{q,\infty}(\Om)$, $q\geq n$, and
initial density is a piecewise constant function. Here we recall
that $L^n(\Om)\subset B^0_{n,\infty}(\Om)$.

\bigskip
Before introducing the main result of the paper we need to make some
conventional notations. Let $[\cdot,\cdot]_\theta$,
$(\cdot,\cdot)_{\theta,r}$ and $(\cdot,\cdot)^0_{\theta,\infty}$ for
$\theta\in (0,1)$, $1\leq r\leq \infty$ be complex, real and
continuous interpolation functors, respectively, see \cite{BL77},
\cite{Tri83} for real and complex interpolation functors, and see
e.g. \cite{Am95}, \S\S 2.4.4, \S2.5, \cite{Tri83} \S\S 1.11.2, page
69 for continuous interpolation functors. We use standard notation
$L^q, H^s_q, B^s_{q,r}$ for Lebesgue spaces, Bessel potential spaces
and Besov spaces, respectively, without distinguishing whether or
not it is the space of scalar-valued functions or vector-valued
functions. For $1<q<\infty$ let
$$L^q_0(\Om):=\{g\in L^q(\Om): \int_{\Om} g\,dx=0\}$$
and let $H^s_{q,0}(\Om)$ for $s>0$ be the completion of
$C^\infty_0(\Om)$ in $H^s_q(\Om)$ and let $b^{s}_{q,\infty}(\Om)$
for $s\in\R$ be the closure of $H^s_q(\Om)$ in
$B^s_{q,\infty}(\Om)$. Let $L^q_\si(\Om)$, $H^{1}_{q,0,\si}(\Om)$,
$b^{0}_{q,\infty,0,\si}(\Om)$, $1<q<\infty$, be the completion of
$C^\infty_{0,\si}(\Om):=\{u\in (C^\infty_{0}(\Om))^n: \div u=0\}$
 in $L^q$-norm, $H^1_q$-norm, $B^0_{q,\infty}$-norm, respectively,
 and
 $$B^{0}_{q,\infty,0,\si}(\Om):=((H^{1}_{q',0,\si}(\Om))', H^{1}_{q,0,\si}(\Om))_{1/2,\infty}.$$
It is known in \cite{RiZh14} that for $1<q<\infty$
$$b^0_{q,\infty,0,\si}(\Om)=\{u\in b^0_{q,\infty}(\Om):\div u=0, u\cdot n|_{\pa\Om}=0\},$$
where the normal component $u\cdot n|_{\pa\Om}$ at the boundary $\pa\Om$
 of $u$ has a meaning in $b^{-1/q}_{q,\infty}(\pa\Om)$,
and $B^0_{q,\infty,0,\si}(\Om)$ is a closed subspace of solenoidal
functions of $B^0_{q,\infty}(\Om)$ and
$$L^q_\si(\Om)\subset b^0_{q,\infty,0,\si}(\Om)\subset  B^0_{q,\infty,0,\si}(\Om),\quad 1<q<\infty.$$
Note that
$$B^0_{q,1}\subset L^q\subset b^0_{q,\infty}\subset  B^0_{q,\infty},\quad 1<q<\infty.$$

Given $\ga\in (0,1]$ and Banach space $X$, we introduce
$L^\infty_\ga(0,T; X):=\{f: t^{1-\ga} f\in L^\infty(0,T;X)\}$ with
norm $\|f\|_{L^\infty_\ga(0,T; X)}:=\|t^{1-\ga}
f(t)\|_{L^\infty(0,T; X)}$. For a set $G$ of $\R^n$ $C_{Lip}(G)$
denotes the set of all Lipschitz continuous functions on $G$ and
$\chi_G$ the characteristic function for $G$. We denote the tensor
product of two tensors $a$, $b$ by $a\otimes b$ and by
$A:B=\sum_{i,j}a_{ij}b_{ij}$ for two matrices $A=(a_{ij})_{1\leq
i,j\leq n}$ and $B=(b_{ij})_{1\leq i,j\leq n}$.

\begin{definition}
 \label{D1.1}
Let $3\leq n\leq q<\infty$. We say that a pair of functions $(\rho,
u)$ is a solution to \eq{1.1} in $(0,T)$ if it satisfies the
followings:
\begin{itemize}
\item[(i)]{
\begin{equation}
\label{E1.8} \rho\in L^\infty(0,T; L^\infty(\Om)),\; u\in
L^\infty(0,T; B^0_{q,\infty}(\Om))\cap L^\infty_{s/2}(0,T;
b^{2-s}_{q,\infty}(\Om))
 \end{equation}
 with some $s\in (0,1)$}, $u|_{\pa\Om}=0$ and $\div u=0$.
\item[(ii)]{Two identities
\begin{equation}
\label{E1.2}
\int_0^T\int_{\Om}(\rho\psi_t+\rho u\cdot \na\psi)\,dxdt
+\int_\Om \rho_0\psi(0,\cdot)\,dx=0,\quad \forall\psi\in C_0^1([0,T)\ti\Om),
\end{equation}
and
\begin{equation}
\begin{array}{l}
\label{E1.3}
 \int_0^T\int_{\Om}[\rho u\cdot\vp_t+\mu u\cdot
\Da\vp+\rho u\otimes u:\na\vp]\,dxdt +\int_\Om
\rho_0u_0\cdot\vp(0,\cdot)\,dx=0,\ek \hfill \forall \vp\in
C_0^\infty([0,T)\ti\Om)^n\; (\div \vp=0).
\end{array}
\end{equation}
hold true.}
\end{itemize}
\end{definition}

If $(\rho, u)$ is a solution to \eq{1.1} in the sense of Definition
\ref{D1.1}, then, by standard argument using De-Rham's lemma, it
follows that there is a distribution $P$ in $\Om$ such that $u$ and
$P$ satisfy the momentum equations of \eq{1.1} in the sense of
distribution and the initial condition in \eq{1.1} is satisfied in a
weak sense. Hence, the triple $(\rho, u, P)$ is also called a
solution to \eq{1.1}.

\bigskip
The main result of the paper is the following:
\begin{theorem}
\label{T1.2} {\rm Let $\Om$ be a bounded domain of $C^2$-class,
$\Om_1$ a subdomain of $\Om$ with Lipschitz boundary and
$\Om_2=\Om\setminus \bar\Om_1$. Let
\begin{equation}
\label{E1.4}
\rho_0(x)=\rho_{01}\chi_{\Om_1}(x)+\rho_{02}\chi_{\Om_2}(x),x\in\Om,
\quad \rho_{02}>\rho_{01}>0,
\end{equation}
and let $u_0\in B^0_{q,\infty,0,\si}(\Om), q\geq n$. Then,
 for any $s\in (0, \min\{\frac{n}{q},1-\frac{1}{q}\})$
 there are some constants $\da_i=\da_i(q,n,s,\Om)>0, i=1,2,$
 such that if $$\frac{\rho_{02}-\rho_{01}}{\rho_{01}}<\da_1,\quad
\|u_0\|_{B^0_{q,\infty}(\Om)}<\frac{\da_2\mu}{\rho_{02}} ,$$
 then  \eq{1.1} has a solution
$(\rho,u, \na P)$ satisfying \eq{1.8}-\eq{1.3} and
\begin{equation}
\label{E1.9}
\begin{array}{l}
u\in L^\infty_{1-\theta+s/2}(0,T;
B^{2\theta-s}_{q,1}(\Om)),\;\forall \theta\in (s/2,1),\quad
  u_t, \na P\in L^\infty_{s/2}(0,T; b^{-s}_{q,\infty}(\Om)).
\end{array}
\end{equation}
 The solution $(\rho,u,\na P)$ is unique if $q>n$ and $s<\frac{q-n}{2q-n}$.
 }
\end{theorem}

\begin{remark}
{\rm (i)  Since one has
$$L^\infty_{s/2}(0,T;
b^{2-s}_{q,\infty}(\Om))\subset L^1(0,T;W^{1,\infty}(\Om)) \subset
L^1(0,T;C_{Lip}(\Om)),\;q>n,$$ for bounded domain $\Om$, it follows
by Theorem \ref{T1.2} that if $(\rho, u)$ is a solution to \eq{1.1}
by Theorem \ref{T1.2} for $q>n$, then the density $\rho(t)$ is
expressed by
$\rho(t,x)=\rho_{01}\chi_{\Om_1(t)}(x)+\rho_{02}\chi_{\Om_2(t)}(x)$,
$t\in (0,T), x\in\Om$, where $\{X(t,\cdot)\}_{t\geq 0}$ stands for
the semiflow transported by the vector field $u$, that is,
$$X(t,y)=y+\int_0^t u(\tau, X(\tau,y))\,d\tau, \quad t\in (0,T), y\in \Om,$$
and $\Om_{i}(t)=\{X(t,y): y\in \Om_{i}\}, i=1,2$. Note that
$X(t,\cdot)$  is  a $C^1$-diffeomorphism over $\Om$ for each $t>0$.
Therefore, the initial interface of two fluids of different density
$\rho_{01}$ and $\rho_{02}$ remains unscattered during the time
period $(0,T)$, which implies that the two fluids are immiscible in
$(0,T)$.

However, if $q=n$, there is no guarantee that the initial interface
may be ``broken" with the evolution of time, and the two fluids may
be mixed together after some time.

(ii) In Theorem \ref{T1.2}, the initial interface of two immiscible
different fluids $\pa\Om_1\cap\pa\Om_2$ belongs to $C^{0,1}$-class
since $\pa\Om_1\in C^{0,1}$ and $\pa\Om_2=\pa\Om -
(\pa\Om\cap\pa\Om_1)$. Indeed, $\pa\Om_1\cap\pa\Om_2$ is allowed to
be a fractal set, for example, a {\it so-called} $d$-set for $d \in
(0,n)$. If this weakened assumption is used in Theorem \ref{T1.2},
the restriction $s\in (0,\min\{\frac{n}{q},1-\frac{n-d}{q}\})$ is
required, see Remark \ref{R4.3} below. }
\end{remark}

In order to prove the main result, beforehand, we consider momentum
equations with a prescribed variable density and transport equations
in Section 2 and 3, respectively. In Section 2, existence of a
solution to instationary Stokes system with nonzero divergence is
proved (Theorem \ref{T2.7}) relying on a consideration of the
divergence problem (Subsection 2.1) and on maximal
$L^\infty_\gamma$-regularity of the Stokes operator in
$b^{\al}_{q,\infty}(\Om)$, $\al\in\R,$ exploited in \cite{RiZh14}.
Then, unique solvability for the nonlinear momentum equations with
prescribed density
\begin{equation}
\label{}
\begin{array}{rcl}
       u_t-\nu\Da u +\na p =a(\nu\Da u -\na p)-(u \cdot \na ) u \,\, &\text{in}&(0,T)\ti\Om,\ek
      \div u = 0 \,\, &\text{in} &(0,T)\ti\Om,\ek
     u= 0\, &\text{on}&(0,T)\ti\pa\Om,\ek
     u(0,x)=u_0 &\text{in}&\Om,
\end{array}
\end{equation}
via standard fixed point argument using linearization when the norms
of $u_0$ and $a\in L^\infty(0,T; {\cal M}(b^{-s}_{q,\infty}(\Om))$
are small (Theorem \ref{T2.4}). Here and in what follows, ${\cal
M}(b^{-s}_{q,\infty}(\Om))$ denotes the space of all pointwise
multipliers for $b^{-s}_{q,\infty}(\Om)$, more precisely,
$$\begin{array}{l}
{\cal M}(b^{-s}_{q,\infty}(\Om)):= \{f: f\vp\in
b^{-s}_{q,\infty}(\Om),\vp\in b^{-s}_{q,\infty}(\Om)\},\ek
\|f\|_{{\cal M}(b^{-s}_{q,\infty}(\Om))} :=\sup_{\vp\in
b^{-s}_{q,\infty}(\Om)}\|f\vp\|_{b^{-s}_{q,\infty}(\Om)}.
\end{array}$$
In Section 3, we prove an existence result for a linear transport
equation in $L^\infty(0,T;L^\infty(\Om))\cap L^\infty(0,T; {\cal
M}(b^{-s}_{q,\infty}(\Om))$ with  piecewise constant initial values
(Proposition \ref{P4.2}).

The proof of Theorem \ref{T1.2} is based on the results of Section 2
and 3.  An iterate scheme for \eq{1.6} is constructed to prove
existence part of Theorem \ref{T1.2} (Section 4.1), while uniqueness
part of Theorem \ref{T1.2} is proved by Lagrangian approach
similarly as in \cite{DaMu12},  \cite{DaZh13} and \cite{HPZ13}, but
we heavily use our result of Theorem \ref{T2.7} on the instationary
Stokes system with nonzero divergence and the result of Lemma
\ref{L4.7} on pointwise multiplication in little Nicolskii spaces
$b^{-s}_{q,\infty}(\Om)$, $s>0$, (Section 4.2).

Throughout the paper, we do not distinguish estimate constants in
the proofs and denote them by $c$ or $C$ as long as no confusion
arises, i.e., the constants $c$ or $C$ may be different from line to
line. We always denote the conjugate number of $q\in (1,\infty)$ by
$q'$, i.e. $q'=q/(q-1)$.

%

\section{Existence for momentum equations}
In this section we consider the existence for the momentum equations
\eq{3.1} with variable density $\rho=\frac{1}{\bar\rho}(1+a)$ fixed.
For the study of the nonlinear problem,  we consider instationary
Stokes problem  with generally a nonzero divergence so as to use the
result of it for the proof of uniqueness of solutions to \eq{1.1},
see \eq{4.35}. To this end, we need to study the divergence problem,
first.
\subsection{Divergence problem} Let $\Om$ be a bounded domain of
$\R^n$, $n\geq 2$, with $C^2$-boundary $\pa\Om$. The divergence
problem
\begin{equation}
\label{E2.7} \div u=g\quad \text{in }\Om, \quad u|_{\pa\Om}=0,
\end{equation}
furnishes important basic tools for the theory of Navier-Stokes
equations and is considerably studied in some references, see e.g.
\cite{Bo79}, \cite{BoSo90}, \cite{FaSo94}-\cite{GHH05}.

In \cite{Bo79}, a solution operator for \eq{2.7},  so-called  {\it
Bogovskii's operator},
  in a star-shaped domain is constructed and, moreover,
  existence of  a solution operator $B$ for \eq{2.7}
  in bounded locally Lipschitz domain satisfying
$$\div Bg=g\quad\text{and}\quad B\in {\cal L}(L^q_0(\Om), H^{1}_{q,0}(\Om)), 1<q<\infty,$$
  is proved.
Moreover, if $\pa\Om$ is smooth enough,
$$B\in {\cal L}(H^{m}_{q,0}(\Om)\cap L^q_0(\Om),
H^{m+1}_{q,0}(\Om)),m\in\N,$$
 (see e.g. \cite{Ga11}, Chapter 3), and, furthermore, if
 $0<s<1-1/q$, then $B$ has a unique extension in  ${\cal L}((H^{s+1}_{q'}(\Om))',(H^{s}_{
q'}(\Om))')$ (cf. \cite{GHH05}).

Note that solutions to \eq{2.7} are not unique. By \cite{FaSo94},
Theorem 1.2 it follows that, given $f\in L^q(\Om)$ and $g\in
H^{1}_{q}(\Om)\cap L^q_0(\Om)$,
 the problem
\begin{equation}
\label{E2.8}
\begin{array}{rcl}
-\Da u+\na p &=& f\quad\text{in }\Om,\ek \div u &=& g\quad\text{in
}\Om,\ek
 u&=& 0\quad\text{on }\pa\Om,\ek
\end{array}
\end{equation}
has a unique solution $u\in H^{2}_{q}(\Om)\cap H^{1}_{q,0}(\Om)$
satisfying
$$\|u\|_{H^{2}_{q}(\Om)\cap H^{1}_{q,0}(\Om)}\leq C(\Om)(\|f\|_{L^q(\Om)}
+\|g\|_{W^{1,q}(\Om)})$$ and, if $f=0$,
$$\|u\|_{L^{q}(\Om)}\leq C(\Om)\|g\|_{(H^{1}_{q'}(\Om)\cap L^{q'}_0(\Om))'}.$$
Therefore, when $f=0$ in \eq{2.8}, the operator ${\cal R}: g\mapsto
u$ defines another solution operator for \eq{2.7}. On the other
hand, it is easily seen by standard argument that
$$\|{\cal R}g\|_{H^1_{q,0}(\Om)}\leq c\|g\|_{L^q_0(\Om)}$$
starting from the existence result for $L^q$-weak solution to
\eq{2.8} with $g=0$, see \cite{Ga11}, Theorem IV. 6.1 (b).
 Thus,
\begin{equation}
\label{E2.20}
\begin{array}{l}
  {\cal R}\in {\cal L}(H^{1}_{q}(\Om)\cap L^q_0(\Om),
H^{2}_{q}(\Om)\cap H^{1}_{q,0}(\Om)),\ek
  {\cal R}\in  {\cal L}(L^q_0(\Om),H^{1}_{q,0}(\Om)),\ek
    {\cal R}\in  {\cal L}((H^{1}_{q'}(\Om)\cap
L^{q'}_0(\Om))',L^q_0(\Om)).
\end{array}
\end{equation}

We shall show that ${\cal R}$ can be continuously extended as an
operator from $(H^{s+1}_{q'}(\Om)\cap L^{q'}_0(\Om))'$ to $(H^{s}_{
q'}(\Om))'$ and for any $s>0$. More precisely, we have

\begin{proposition}
\label{L2.1} {\rm Let $s>0$ and $\Om\subset\R^n, n\geq 2,$ be a
bounded domain with $C^{2+s}$-boundary and $q\in (1,\infty)$. Then
the solution operator ${\cal R}$ for \eq{2.7} can be uniquely
extended as
\begin{equation}
\label{E2.15} {\cal R}\in {\cal L}((H^{s+1}_{ q'}(\Om)\cap
L^{q'}_0(\Om))',(H^{s}_{ q'}(\Om))').
\end{equation}
  Moreover, if $\pa\Om\in C^{3+s}$, then
\begin{equation}
\label{E2.16}
  {\cal R}\in {\rm Isom}((H^{s+1}_{ q'}(\Om)\cap L^{q'}_0(\Om))',(H^{s}_{ q'}(\Om))').
  \end{equation}
   }
\end{proposition}
{\bf Proof:} The proof relies on a duality argument based on the
regularity for the Stokes system \eq{2.8}. Consider the problem
\begin{equation}
\label{E2.12}
\begin{array}{rcl}
-\Da z+\na \zeta &=& \vp\quad\text{in }\Om,\ek \div z &=&
0\quad\text{in }\Om,\ek
 z&=& 0\quad\text{on }\pa\Om.
\end{array}
\end{equation}
By the well-known regularity theory for stationary Stokes systems,
see \cite{Ga11}, for any $\vp\in H^{s}_{q'}(\Om)$ \eq{2.12} has a
unique solution $\{z,\zeta\}\in (H^{2+s}_{q'}(\Om)\cap
H^{1}_{q',0}(\Om))\ti (H^{1+s}_{q'}(\Om)\cap L^q_0(\Om))$ satisfying
$$\|z\|_{H^{2+s}_{q'}(\Om)}+\|\zeta\|_{H^{1+s}_{q'}(\Om)}\leq c \|\vp\|_{H^{s}_{q'}(\Om)}.$$
Now, define the operator
\begin{equation}
\label{E2.17} S\in {\cal L}(H^{s}_{q'}(\Om), H^{1+s}_{q'}(\Om)\cap
L^{q'}_0(\Om)),\quad S\vp:=\zeta.
\end{equation}
Then we get for solution $u\in H^{1}_{q,0}(\Om)$ to \eq{2.8} with
$f=0$, $g\in L^q_0(\Om)$ that
\begin{equation}
\label{E2.13}
\begin{array}{rcl}
\lan u, \vp\ran_{(H^{s}_{q'}(\Om))',H^{s}_{q'}(\Om)}
 &=&\lan u, \vp\ran_{L^q(\Om),L^{q'}(\Om)}
 =\lan u, -\Da z+\na\zeta\ran_{L^q(\Om),L^{q'}(\Om)}\ek
 &=&\lan -\Da u+\na p, z\ran_{H^{-1}_{q}(\Om),H^{1}_{q',0}(\Om)}
            -\lan \div u, \zeta\ran_{L^q(\Om),L^{q'}(\Om)}\ek
 &=&\lan f, z\ran_{H^{-1}_{q}(\Om),H^{1}_{q',0}(\Om)}-\lan g, \zeta\ran_{(H^{s+1}_{q'}(\Om))',H^{s+1}_{q'}(\Om)}\ek
 &=&-\lan g,S\vp\ran_{(H^{s+1}_{q'}(\Om)\cap L^{q'}_0(\Om))',H^{s+1}_{q'}(\Om)\cap L^{q'}_0(\Om)}, \forall \vp\in
H^{s}_{q'}(\Om);
\end{array}
\end{equation}
here note that
$$(H^{s+1}_{q'}(\Om)\cap L^{q'}_0(\Om))'=\{g\in (H^{s+1}_{q'}(\Om))': \lan g,1\ran_\Om=0\}$$
and
\begin{equation}
\label{E2.14} H^{s_1}_{q}\Om\cap L^{q}_0(\Om)\dense
L^q_0(\Om)=(L^{q'}_0(\Om))'\dense (H^{s_2}_{q'}\Om\cap
L^{q'}_0(\Om))'
\end{equation}
for any $s_1, s_2>0$. From \eq{2.13} we have $u=S'g$, where
$$S'\in {\cal L}((H^{1+s}_{q'}(\Om)\cap L^{q'}_0(\Om))', (H^{s}_{q'}(\Om))')$$
is the dual operator of $S$. Thus we have
  $${\cal R}=S'|_{L^{q}_0(\Om)}.$$
 It is clear  in view of \eq{2.14} that $S'$
is the unique extension of ${\cal R}$, in other words, ${\cal R}$
can be uniquely extended as \eq{2.15}.

Now, assume that $\pa\Om\in C^{3+s}$. By the well-known uniqueness
of regular solution to \eq{2.12} (cf. \cite{Ga11}, \cite{Te77}), the
operator $S$ in \eq{2.17} is injective. Therefore, if we prove that
$S$ is surjective, then \eq{2.16} is proved. Given arbitrary $y\in
H^{s+1}_{q'}(\Om)\cap L^{q'}_0(\Om)$, let $\eta_0\in
H^{s+1}_{q'}(\Om)$ be the (unique) solution to the problem
$$-\Da\eta_0=\Da y\,\,\text{in } \Om,\quad \eta_0|_{\pa\Om}=0,$$
and let $z_0\in H^{2+s}_{q'}(\Om)\cap H^{1}_{q',0}(\Om)$ be the
(unique) solution to
$$-\Da z_0=\na(\eta_0- y)\,\,\text{in } \Om,\quad z_0|_{\pa\Om}=0.$$
Note that $\div z_0\in H^{1+s}_{q'}(\Om)\cap L^{q'}_0(\Om)$ and
$\div z_0\neq 0$, in general. Let $z_1=\na\eta_1$, where $\eta_1\in
H^{3+s}_{q'}(\Om)$ is the solution to
$$-\Da \eta_1=-\div z_0\,\,\text{in } \Om,\quad \eta_1|_{\pa\Om}=0,$$
and $z_2\in H^{2+s}_{q'}(\Om)$ the solution to
$$-\Da z_2+\na \pi =0\,\,\text{in } \Om,\quad \div z_2=0,\quad z_2|_{\pa\Om}=-(z_0+z_1)|_{\pa\Om}.$$
Then, $z:=z_0+z_1+z_2$ and $y$ solve the system \eq{2.12} with
$\vp=\na\eta_0-\Da z_1-\Da z_2\in H^{s}_{q'}(\Om)$, that is,
 $S \vp=y$.

The proof of the lemma is complete. \qed 
\begin{remark}
\label{R2.2} {\rm By Proposition \ref{L2.1}, \eq{2.16} and its
proof, we can conclude that, if $\Om$ is a bounded domain with
sufficiently smooth boundary, then the solution operator ${\cal R}$
for the divergence problem is expressed by
 ${\cal R}= Q'\in {\rm Isom}((H^{s+1}_{ q'}(\Om)\cap L^{q'}_0(\Om))',(H^{s}_{ q'}(\Om))')$,
 where $Q'$ is dual of the linear operator $Q$ defined by \eq{2.17}. Moreover,
$${\cal R}^{-1}\in {\rm Isom}((H^{s}_{ q'}(\Om))', (H^{s+1}_{ q'}(\Om)\cap L^{q'}_0(\Om))').$$
Note that if $s>1-1/q$, then $(H^{s}_{ q'}(\Om))'$ is not included
in the space of Schwartz distributions. Hence, ${\cal R}^{-1}u$ for
$u\in (H^{s}_{ q'}(\Om))'$, $s>1-1/q$, may not be regarded as $\div
u$ in the sense of distribution.
 }
\end{remark}

\begin{proposition}
\label{C2.3} {\rm Let $\Om$ be a bounded domain  of $\R^n$, $n\geq
2$, with $C^2$-boundary, $q\in (1,\infty)$ and let ${\cal R}$ be the
same as in Proposition \ref{L2.1}. Then,
\begin{equation}
\label{E2.9} {\cal R}\in {\cal L}(H^{\al}_{q}(\Om)\cap L^q_0(\Om),
H^{1+\al}_{q}(\Om)
  \cap H^{1}_{q,0}(\Om)), \al\in [0,1],
\end{equation}
\begin{equation}
\label{E2.10} {\cal R}\in {\cal L}(B^{\al}_{q,r}(\Om)\cap
L^q_0(\Om), B^{1+\al}_{q,r}(\Om)\cap H^{1}_{q,0}(\Om)), \al\in
(0,1), 1\leq r\leq\infty,
\end{equation}
and
\begin{equation}
\label{E2.22} {\cal R}\in {\cal L}(b^{\al}_{q,\infty}(\Om)\cap
 L^q_0(\Om), b^{1+\al}_{q,\infty}(\Om)\cap H^{1}_{q,0}(\Om)),\al\in (0,1).
\end{equation}

Moreover, if $\pa\Om\in C^{2+s}$, $s>0$, then one has, in addition
to \eq{2.15},
\begin{equation}
\label{E2.19} {\cal R}\in {\cal L}((B^{s+1}_{ q',r}(\Om)\cap
L^{q'}_0(\Om))',(B^{s}_{ q',r}(\Om))') \cap {\cal L}((b^{s+1}_{
q',\infty}(\Om)\cap L^{q'}_0(\Om))',(b^{s}_{ q',\infty}(\Om))')
\end{equation}
 for any $r\in [1,\infty)$.
 }
\end{proposition}
{\bf Proof:} Note that  \eq{2.9} for the case $\al=0$ or $\al=1$ has
already been proved by \eq{2.20}.  We get by \eq{2.20} that
$${\cal R}\in
{\cal L}((L^q_0(\Om), H^1_q(\Om)\cap L^q_0(\Om))_{\al},
(H^1_q(\Om),H^2_q(\Om))_{\al}), \forall \al\in (0,1),
$$
 where $(\cdot,\cdot)_{\al}$ is any interpolation functor.
Since $L^q_0(\Om)$ is a subspace of $L^q(\Om)$, which is the range
of the projection operator $\Pi$ defined on $L^q(\Om)$ by $\Pi
g:=g-\int_{\Om}g\,dx$, one gets by interpolation theory, see e.g.
\cite{Tri83}, Theorem 1.17.1/1, that
$$(L^q_0(\Om), H^1_q(\Om)\cap L^q_0(\Om))_{\al}=
 (L^q(\Om), H^1_q(\Om))_{\al}\cap L^q_0(\Om), \al\in (0,1).$$
Then, due to the relations
 $$\begin{array}{l} [H^{s_1}_q(\Om),
H^{s_2}_q(\Om)]_{\al}=H^{(1-\al)s_1+\al s_2}_q(\Om),\ek
 (H^{s_1}_q(\Om),
H^{s_2}_q(\Om))_{\al,r}=B^{(1-\al)s_1+\al s_2}_{q,r}(\Om),\; 1\leq
r\leq \infty, \end{array}$$
 and
$$(H^{s_1}_q(\Om), H^{s_2}_q(\Om))^0_{\al,\infty}=b^{(1-\al)s_1+\al
s_2}_{q,\infty}(\Om),$$
  we get \eq{2.9}, \eq{2.10} and \eq{2.22} by
taking $(\cdot,\cdot)_{\al}$ as complex, real and continuous
interpolation functors, respectively.

Finally, let us prove \eq{2.19}. Due to \eq{2.15}, one gets by
interpolation that
\begin{equation}
\label{E2.11} {\cal R}\in {\cal L}(((H^{s_1+1}_{ q'}(\Om)\cap
L^{q'}_0(\Om))',(H^{s_2+1}_{ q'}(\Om) \cap
L^{q'}_0(\Om))')_\theta,\, ((H^{s_1}_{ q'}(\Om))',(H^{s_2}_{
q'}(\Om))')_\theta), \theta\in (0,1),
\end{equation}
 where $(\cdot,\cdot)_\theta$
is any interpolation functor on  couple of Banach space. Recall the
relation of dual interpolation
$$\big((X_1,X_2)_{\theta,r}\big)'=(X'_1,X'_2)_{\theta,r'},\quad 1\leq
r<\infty,$$
 and
$$\big((X_1,X_2)^0_{\theta,\infty}\big)'=(X'_1,X'_2)_{\theta,1}$$
(see \cite{BL77}, Theorem 3.7.1 or \cite{Tri83}, Theorem 1.11.2, or
\cite{Am95}, \S2.6). Hence,
$$\begin{array}{l}
((H^{s_1+1}_{ q'}(\Om)\cap L^{q'}_0(\Om))',(H^{s_2+1}_{ q'}(\Om)\cap
L^{q'}_0(\Om))')_{\theta,r'}\ek \qquad
 = \big((H^{s_1+1}_{q'}(\Om)\cap L^{q'}_0(\Om),H^{s_2+1}_{ q'}(\Om)\cap
L^{q'}_0(\Om))_{\theta,r}\big)' \ek \qquad
 =((H^{s_1+1}_{q'}(\Om),H^{s_2+1}_{q'}(\Om))_{\theta,r}\cap L^{q'}_0(\Om))'
\end{array}$$
and, similarly,
$$((H^{s_1+1}_{ q'}(\Om)\cap L^{q'}_0(\Om))',(H^{s_2+1}_{ q'}(\Om)\cap
L^{q'}_0(\Om))')_{\theta,1} =((H^{s_1+1}_{ q'}(\Om),H^{s_2+1}_{
q'}(\Om))^0_{\theta,\infty}\cap L^{q'}_0(\Om))'.
$$
 Thus, taking $(\cdot,\cdot)_\theta$ in \eq{2.11} as
$(\cdot,\cdot)_{\theta,r'}, 1\leq r< \infty,$ or
$(\cdot,\cdot)_{\theta, 1}$, we get \eq{2.19}. \qed
\begin{proposition}
\label{R2.4} {\rm For $s\in (0,1)$ and $q\in (1, \infty)$ one has
\begin{equation}
\label{E2.18n} \|{\cal R}\div h\|_{b^{-s}_{q,\infty}(\Om)} \leq
c(\Om,q)\|h\|_{b^{-s}_{q,\infty}(\Om)},\forall h\in
b^{-s}_{q,\infty}(\Om).
\end{equation}
}
\end{proposition}
{\bf Proof:} For any  $h\in H^1_{q,0}(\Om)$ and $\psi\in
H^{s+1}_{q'}(\Om)\cap L^{q'}_0(\Om)$ one has
$$|\lan \div h,\psi\ran_\Om|=|\lan h,\na\psi\ran_\Om|
\leq \|h\|_{(H^s_{q'}(\Om))'}\|\psi\|_{H^{s+1}_{q'}(\Om)\cap
L^{q'}_0(\Om)},$$ which yields
$$\|\div h\|_{(H^{s+1}_{q'}(\Om)\cap L^{q'}_0(\Om))'}\leq \|h\|_{(H^{s}_{q'}(\Om))'},
                  \quad\forall h\in (H^{s}_{q'}(\Om))',$$
in view of the denseness of $H^1_{q,0}(\Om)$ in
$(H^{s}_{q'}(\Om))'$. Therefore, by Proposition \ref{L2.1} we have
\begin{equation}
\begin{array}{rcl}
\label{E2.18} \|{\cal R}\div h\|_{(H^{s}_{q'}(\Om))'} &\leq &
c(\Om,q)\|\div h\|_{(H^{s+1}_{q'}(\Om)\cap L^{q'}_0(\Om))'}\ek
  &\leq & c(\Om,q)\|h\|_{(H^{s}_{q'}(\Om))'},\;\forall h\in
(H^{s}_{q'}(\Om))',
\end{array}
\end{equation}
with $c=c(\Om,q)$. Using the duality theorem for real interpolation
(e.g. \cite{BL77}, Theorem 3.7.1), one gets for
$s=(1-\theta)s_1+\theta s_2$, $s_1, s_2\in (0,1), \theta\in (0,1)$,
that
$$((H^{s_1}_{q'}(\Om))', (H^{s_2}_{q'}(\Om))')_{\theta,\infty}
 =(H^{s_1}_{q'}(\Om), H^{s_2}_{q'}(\Om))'_{\theta,1}=(B^{s}_{q',1}(\Om))'=b^{-s}_{q,\infty}(\Om).$$
Thus,  \eq{2.18n} follows from \eq{2.18} by real interpolation
$(\cdot,\cdot)_{\theta,\infty}$. \qed

\subsection{Instationary Stokes problem}

\par\bigskip
Consider the initial boundary value problem
for instationary Stokes system
\begin{equation}
\label{E2.21}
\begin{array}{rcl}
     u_t-\nu\Da u  +\na p  = f \,\, &\text{in}&(0,T)\ti\Om,\ek
      \div u = g \,\, &\text{in} &(0,T)\ti\Om,\ek
     u= 0\, &\text{on}&(0,T)\ti\pa\Om,\ek
     u(0,x)=u_0 &\text{in}&\Om,
\end{array}
\end{equation}
where $\Om$ is a bounded domain of $\R^n, n\geq 3$. 

\begin{lemma}
\label{L2.1n} {\rm Let $\Om$ be a bounded domain of $\R^n, n\geq 2,$
with $C^1$-boundary. Then,
\begin{equation}
\label{E2.3n}
 \|g\|_{b^{1-s}_{q,\infty}(\Om)}\leq c(\Om,q)\|\na
g\|_{b^{-s}_{q,\infty}(\Om)},\quad \forall g\in
b^{1-s}_{q,\infty}(\Om)\cap L^q_0(\Om).
 \end{equation}
}
\end{lemma}
{\bf Proof:} By Poincar\'e's inequality, one gets
\begin{equation}
\label{E2.1n}
 \|g\|_{H^1_q(\Om)}\leq c(\Om,q)\|\na g\|_{L^q(\Om)}\leq c(\Om,q)\|g\|_{H^1_q(\Om)},\quad \forall g\in
H^{1}_{q}(\Om)\cap L^q_0(\Om).
 \end{equation}

  On the other hand, one gets
\begin{equation}
\label{E2.2n} \|g\|_{L^q_0(\Om)}\leq c(\Om,q)\|\na
g\|_{H^{-1}_{q}(\Om)}\leq  c(\Om,q)\|g\|_{L^q_0(\Om)}.
 \end{equation}
In fact, for any $\psi\in L^{q'}_0(\Om)$ there is some $h\in
H^1_{q',0}(\Om)$ such that
$$\div h=\psi,\quad \|\na h\|_{q'}\leq c(q,\Om)\|\psi\|_{q'}$$
(see e.g. \cite{Ga11}, Theorem III.3.1). Therefore,
$$|\langle g, \psi\rangle|=|\langle \na g, h\rangle|\leq \|\na g\|_{H^{-1}_{q}(\Om)}\|h\|_{H^{1}_{q',0}(\Om)}\leq
c(q,\Om)\|\na g\|_{H^{-1}_{q}(\Om)}\|\psi\|_{q'},$$ which yields
\eq{2.2n}.

Thus, using continuous interpolation
$$(L^q(\Om),
H^1_q(\Om))^0_{1-s,\infty}=b^{1-s}_{q,\infty}(\Om),\quad
(H^{-1}_q(\Om), L^q(\Om))^0_{1-s,\infty}=b^{-s}_{q,\infty}(\Om),$$
the assertion of the lemma follows from \eq{2.1n} and \eq{2.2n}.\qed

\begin{theorem}
\label{T2.7} {\rm Let $\Om\subset \R^n, n\geq 3,$ be a bounded
domain of $C^2$-class and $u_0\in B^{0}_{q,\infty,0,\si}(\Om)$. Let
$f\in L^\infty_{s/2}(0,T;b^{-s}_{q,\infty}(\Om))$, $g\in
L^\infty_{s/2}(0,T;b^{1-s}_{q,\infty}(\Om)\cap L^q_0(\Om))$, $q\in
(1,\infty)$,
 $s\in (0,1)$,
   and $g=\div R$ with some distribution $R=R(t,\cdot)$
  such that $R_t\in L^\infty_{s/2}(0,T;b^{-s}_{q,\infty}(\Om))$, $R(0)=0$.
Then the problem \eq{2.21} has a unique solution $\{u,\na p\}$ such
that
\begin{equation}
\label{E2.24}
\begin{array}{l}
u\in L^\infty(0,T; B^{0}_{q,\infty}(\Om))\cap L^\infty_{s/2}(0,T;
b^{2-s}_{q,\infty}(\Om)),\ek
 u_t, \na p\in L^\infty_{s/2}(0,T; b^{-s}_{q,\infty}(\Om)),\ek
  \|u_t, \nu \na^2 u, \na p\|_{L^\infty_{s/2}(0,T; b^{-s}_{q,\infty}(\Om))}
  +\nu^{s/2}\|u\|_{L^\infty(0,T; B^{0}_{q,\infty}(\Om))}
  \ek
\hspace{4cm}    \leq c(\|f, \nu\na g,R_t\|_{L^\infty_{s/2}(0,T;b^{-s}_{q,\infty}(\Om))}
  +\nu^{s/2}\|u_0\|_{B^{0}_{q,\infty}(\Om)})
\end{array}
\end{equation}
with constant $c>0$ depending only on $q,n,s,\Om$  and independent of $T$.

Moreover, the solution $u$ to \eq{2.21} satisfies
 \begin{equation}
\label{E2.31}
\begin{array}{l}
u\in L^\infty_{1-\theta+s/2}(0,T; B^{2\theta-s}_{q,1}(\Om)),\ek
\nu^\theta\|u\|_{L^\infty_{1-\theta+s/2}(0,T;
B^{2\theta-s}_{q,1}(\Om))}
  \leq
  c(\|f, \nu\na g,R_t\|_{L^\infty_{s/2}(0,T;b^{-s}_{q,\infty}(\Om))}
    +\nu^{s/2}\|u_0\|_{B^{0}_{q,\infty}(\Om)}), \ek
    \hfill \forall \theta\in (s/2,1),\end{array}
\end{equation}
with constant $c>0$ depending only on $q,n,s,\Om$ and independent of
$T$ and $\nu$.

}

\end{theorem}

{\bf Proof:} First let us prove the theorem for $\nu=1$.

 The case $g=0$ has already been proved by \cite{RiZh14}, Corollary
 4.14, (i) as a particular case.

 Let $g\in L^{\infty}_{s/2}(0,T; b^{1-s}_{q,\infty}(\Om))$
 be identically not $0$ and  $w(t)={\cal R}g(t), t\in (0,T)$,
where ${\cal R}$ is the solution operator for \eq{2.7} constructed
through the solution to \eq{2.8} with $f=0$. Then, $\div w(t)=g(t)$
for all $t\in (0,T)$ and it is easily checked by Proposition
\ref{C2.3}, \eq{2.22} that
\begin{equation}
\label{E2.4n}
\begin{array}{l}
 w \in L^\infty_{s/2}(0,T;b^{2-s}_{q,\infty}(\Om)\cap
H^1_{q,0}(\Om)),\ek
 \|w\|_{L^\infty_{s/2}(0,T; b^{2-s}_{q,\infty}(\Om))} \leq
c\|g\|_{L^\infty_{s/2}(0,T;b^{1-s}_{q,\infty}(\Om))}.
\end{array}
\end{equation}

Moreover, since $w_t={\cal R}g_t={\cal R}\div R_t$ and $R_t\in
L^\infty_{s/2}(0,T;b^{-s}_{q,\infty}(\Om))$, we get by Proposition
\ref{R2.4}, \eq{2.18n} that
\begin{equation}
\label{E2.7n} w_t\in L^\infty_{s/2}(0,T;b^{-s}_{q,\infty}(\Om)), \ek
\|w_t\|_{L^\infty_{s/2}(0,T;b^{-s}_{q,\infty}(\Om))}\leq
c\|R_t\|_{L^\infty_{s/2}(0,T;b^{-s}_{q,\infty}(\Om))}
\end{equation}

Furthermore, since $R(0)=0$ and $R_t\in L^\infty_{s/2}(0,T;
b^{-s}_{q,\infty}(\Om))$, we have $R(t)=\int_0^t
R_{\tau}\,d\tau=\int_0^t
\tau^{-1+s/2}(\tau^{1-s/2}R_{\tau})\,d\tau$, and
 $$\|R(t)\|_{b^{-s}_{q,\infty}(\Om)}\leq \frac{2t^{s/2}}{s}\|R_t\|_{
L^\infty_{s/2}(0,T; b^{-s}_{q,\infty}(\Om))},\;\forall t\in (0,T).$$
 Hence,
 $$\begin{array}{l}
 t^{-s/2} R\in L^\infty(0,T;
b^{-s}_{q,\infty}(\Om)),\ek
 \|t^{-s/2}R\|_{L^\infty(0,T;
b^{-s}_{q,\infty}(\Om))}\leq \frac{2}{s}\|R_t\|_{
L^\infty_{s/2}(0,T; b^{-s}_{q,\infty}(\Om))},
 \end{array}$$
which yields by Proposition \ref{R2.4}, \eq{2.18n} that
\begin{equation}
\label{E2.5n}
\begin{array}{l}
 t^{-s/2}w=t^{-s/2}{\cal R}\div R\in L^\infty(0,T;
b^{-s}_{q,\infty}(\Om)),\ek
  \|t^{-s/2}w\|_{L^\infty(0,T;
b^{-s}_{q,\infty}(\Om))}\leq \frac{2}{s}\|R_t\|_{
L^\infty_{s/2}(0,T; b^{-s}_{q,\infty}(\Om))}.
\end{array}
\end{equation}
Therefore, it follows from \eq{2.4n}, \eq{2.5n} that, by complex
interpolation,
 $$ \begin{array}{l}
 w \in  L^\infty_{1-s/2}(0,T;
[b^{2-s}_{q,\infty}(\Om), b^{-s}_{q,\infty}(\Om)]_{(2-s)/2})\subset
L^\infty_{1-s/2}(0,T; B^0_{q,\infty}(\Om)),\ek
 \|w(t)\|_{B^{0}_{q,\infty}(\Om)}\leq
c(s)\|w(t)\|^{s/2}_{b^{2-s}_{q,\infty}(\Om)}
 \|w(t)\|^{1-s/2}_{b^{-s}_{q,\infty}(\Om)}\ek
 \leq c(s)(t^{(-1+s/2)}\|g\|_{L^\infty_{s/2}(0,T;b^{1-s}_{q,\infty}(\Om)})^{s/2}
 (t^{s/2}\|R_t\|_{L^\infty_{1-s/2}(0,T;b^{-s}_{q,\infty}(\Om))})^{1-s/2}\ek
 =
c(s)\|g\|_{L^\infty_{s/2}(0,T;b^{1-s}_{q,\infty}(\Om))}^{s/2}
 \|R_t\|_{L^\infty_{1-s/2}(0,T;b^{-s}_{q,\infty}(\Om))}^{1-s/2},
 \quad \text{for a.a. } t\in (0,T),
\end{array}
$$
where we used the interpolation relation $[b^{2-s}_{q,\infty}(\Om),
b^{-s}_{q,\infty}(\Om)]_{(2-s)/2}\subset B^0_{q,\infty}(\Om)$ in
view of $(2-s)(1-\frac{2-s}{2})-s\cdot\frac{2-s}{2}=0$. Thus,
\begin{equation}
\label{E2.6n}
\begin{array}{l}
w\in L^\infty(0,T;B^{0}_{q,\infty}(\Om)),\ek
\|w\|_{L^\infty(0,T;B^{0}_{q,\infty}(\Om))}\leq
c(s)\|g\|_{L^\infty_{s/2}(0,T;b^{1-s}_{q,\infty}(\Om))}^{s/2}
 \|R_t\|_{L^\infty_{1-s/2}(0,T;b^{-s}_{q,\infty}(\Om))}^{1-s/2}.
\end{array}
\end{equation}


Now, introducing the new unknown $U=u-w$, the problem
\eq{2.21} is reduced to a divergence-free problem, that is,
$$\begin{array}{rlccl}
     U_t-\Da U  +\na p  &= &F \,\, &\text{in}&(0,T)\ti\Om,\ek
      \div U &=& 0 \,\, &\text{in} &(0,T)\ti\Om,\ek
     U &=& 0 \,\, &\text{on} &(0,T)\ti\pa\Om,\ek
    U(0,x)&=&U_0&\text{in}&\Om,
\end{array}$$
where $F:=f-w_t-\Da w\in L^\infty_{s/2}(0,T;
b^{-s}_{q,\infty}(\Om))$, $U_0:=u_0-w(0)=u_0\in
B^{0}_{q,\infty,\si}(\Om)$. Thus, by the assertion of the theorem
for the divergence-free case, we get that
$$\begin{array}{l}
U_t, \na^2 U, \na p\in  L^\infty_{s/2}(0,T; b^{-s}_{q,\infty}(\Om)),\,\,U|_{\pa\Om}=0,\quad
U\in L^\infty(0,T; B^{0}_{q,\infty}(\Om)),\ek
  \|U_t,\na^2 U, \na p\|_{L^\infty_{s/2}(0,T; b^{-s}_{q,\infty}(\Om))}
  +\|U\|_{L^\infty(0,T; B^{0}_{q,\infty}(\Om))}
  \ek
\hspace{4cm}    \leq c(\|F\|_{L^\infty_{s/2}(0,T;b^{-s}_{q,\infty}(\Om))}
  +\|U_0\|_{B^{0}_{q,\infty}(\Om)}).
\end{array}
$$
Thus, in view of \eq{2.4n}, \eq{2.7n}, \eq{2.6n} and Lemma
\ref{L2.1n}, we get \eq{2.24} for $u=U+w$.

The uniqueness of solution to \eq{2.21} is clear from \eq{2.24}.

The estimate \eq{2.31} for $\nu=1$ follows from the estimate of $u$
in $L^\infty_{s/2}(0,T; b^{2-s}_{q,\infty}(\Om))\cap L^\infty(0,T;
B^{0}_{q,\infty}(\Om))$ (see \eq{2.24})  by real interpolation using
$(B^{s_1}_{q,\infty}(\Om),
B^{s_1}_{q,\infty}(\Om))_{\gamma,1}=B^{(1-\gamma)s_1+\gamma
s_2}_{q,1}(\Om)$ for all $\gamma \in (0,1)$.

Finally, assuming $\nu\neq 1$, let us prove \eq{2.24} and \eq{2.31}.

Notice that the rescaling transform
\begin{equation}
\label{E2.33}
  (\tilde{u},\tilde{g})(t,x):=\nu^{-1} (u,g)( \nu^{-1}t,x),\,\,
   (\tilde{p}, \tilde{f})(t,x):=\nu^{-2} (p,f)(\nu^{-1}t,x),\,\,
 \tilde{u}_0(x):=\nu^{-1} u_0(x)
\end{equation}
reduces \eq{2.21} with $\nu\neq 1$ to the case with $\nu=1$ and that
$\tilde{g}_t(t,x)=\div \tilde{R}(t,x)$, where
 $\tilde{R}(t,x):=\frac{1}{\nu^2} R(\frac{t}{\nu},x)$.
Then we have
$$\begin{array}{l}
  \|\tilde{u}_t, \na^2 \tilde{u}, \na \tilde{p}\|_{L^\infty_{s/2}(0,\nu T; b^{-s}_{q,\infty}(\Om))}
  +\|\tilde{u}\|_{L^\infty(0,\nu T; B^{0}_{q,\infty}(\Om))}\ek
\hspace{4cm}    \leq c(\|\tilde{f},\na \tilde{g},
\tilde{R}_t\|_{L^\infty_{s/2}(0,\nu T;b^{-s}_{q,\infty}(\Om))}
  +\|\tilde{u}_0\|_{B^{0}_{q,\infty}(\Om)})
\end{array}
$$
with $c=c(q,n,s)$ independent of $\nu$ and $T$. Here we have
$$\begin{array}{l}
\|\tilde{u}_t, \na^2 \tilde{u},\na\tilde{p}\|_{L^\infty_{s/2}(0,\nu T; b^{-s}_{q,\infty}(\Om))}\ek
=
\frac{1}{\nu^2}\nu^{1-s/2}
\|(\frac{t}{\nu})^{1-s/2}u_t(\frac{t}{\nu},x), (\frac{t}{\nu})^{1-s/2}\nu\na^2u(\frac{t}{\nu},x),
(\frac{t}{\nu})^{1-s/2}\na p(\frac{t}{\nu},x)\|_{L^\infty(0,\nu T; b^{-s}_{q,\infty}(\Om))}\ek
=\nu^{-1-s/2}
\|u_t, \nu\na^2u,\na p\|_{L^\infty_{s/2}(0,T; b^{-s}_{q,\infty}(\Om))},
\end{array}$$
$$\|\tilde{u}\|_{L^\infty(0,\nu T; B^{0}_{q,\infty}(\Om))}=
\nu^{-1}\|u\|_{L^\infty(0,T; B^{0}_{q,\infty}(\Om))},
\|\tilde{u}_0\|_{B^{0}_{q,\infty}(\Om)}=\nu^{-1}\|u_0\|_{B^{0}_{q,\infty}(\Om)}
$$
and
  $$\|\tilde{f},\na \tilde{g}, \tilde{R}_t\|_{L^\infty_{s/2}(0,\nu T;b^{-s}_{q,\infty}(\Om))}
  =\nu^{-1-s/2}\|f,\nu\na g, R_t\|_{L^\infty_{s/2}(0,T;b^{-s}_{q,\infty}(\Om))}.
  $$
  Hence \eq{2.24} is proved. In the same way, \eq{2.31} can be proved, and we omit the proof.
\qed

\subsection{Momentum equations with  fixed variable density}
\par
Let $\rho_{-}\leq \rho\leq \rho_{+}$ with some positive constants
$\rho_{-},\rho_{+}$ and let $\frac{1}{\rho}=\frac{1}{\bar\rho}(1+a)$
with constant $\bar\rho\in [\rho_{-}, \rho_{+}]$ and $a=a(t,x)\in
L^\infty((0,T)\ti\Om)$. Then the momentum equation in \eq{1.1} may
be formally written as
\begin{equation}
\label{E3.1}
\begin{array}{rcccl}
     u_t-\nu\Da u  + (u \cdot \na) u+\na p  &=& a(\nu\Da u-\na p) \,\, &\text{in}&(0,T)\ti\Om,\ek
      \div u &=& 0 \,\, &\text{in} &(0,T)\ti\Om,\ek
     u &=& 0\, &\text{on}&(0,T)\ti\pa\Om,\ek
  u(0,x) &=& u_0 &\text{in}&\Om,
\end{array}
\end{equation}
where $0<T\leq \infty$.
\par
\bigskip
We need the following lemma for the proof of Theorem \ref{T2.4},
which is the main result of this section.
\begin{lemma} (see \cite{RiZh14}, Lemma 5.1)
\label{L3.0} {\rm Let $\Om$ be a bounded domain of $\R^n,n\geq 3,$
of $C^2$-class, then for $q\in [n,\infty)$ and $s\in
(0,\frac{n}{q})$
\begin{equation}
\label{E4.4n}
H^{(2-s)(1-\theta)}_{q}(\Om)\cdot H^{(2-s)\theta-1}_{q}(\Om)
\hookrightarrow H^{-s}_{q}(\Om),
\quad \forall \theta\in (\frac{1}{2-s},\frac{1-s+n/q}{2-s}).
\end{equation}
 }
\end{lemma}
\begin{theorem}
\label{T2.4} {\rm Let $\Om\subset\R^n,n\geq 3,$ be a bounded domain
of $C^2$-class. Let $q\geq n$, $a\in L^\infty(0,T;{\cal
M}(b^{-s}_{q,\infty}(\Om)))$ for some $s\in (0,\frac{n}{q})$, and
let $u_0\in B^0_{q,\infty,0,\si}(\Om)$.

Then there are some constants $\da_i=\da_i(q,n,\Om,s)>0, i=1,2,$ and
$M=M(q,n, s,\Om)>0$ independent of $T$ such that if
$$\|a\|_{L^\infty(0,T;{\cal M}(b^{-s}_{q,\infty}(\Om)))}<\da_1,\quad
\|u_0\|_{B^0_{q,\infty}(\Om)}<\da_2\nu,$$
 then  \eq{3.1} has a solution
$(u,\na p)$ satisfying
$$\begin{array}{l}
u\in L^\infty(0,T; B^{0}_{q,\infty}(\Om))\cap
L^\infty_{1-\theta+s/2}(0,T; B^{2\theta-s}_{q,1}(\Om))\cap
L^\infty_{s/2}(0,T; b^{2-s}_{q,\infty}(\Om)),\ek
 u_t, (u\cdot\na)u, \na p\in L^\infty_{s/2}(0,T; b^{-s}_{q,\infty}(\Om))
\end{array}$$
and
\begin{equation}
\label{E3.6}
\begin{array}{l}
\nu^{s/2}\|u\|_{L^\infty(0,T; B^{0}_{q,\infty}(\Om))}
+\nu^\theta\|u\|_{L^\infty_{1-\theta+s/2}(0,T;
B^{2\theta-s}_{q,1}(\Om))}\ek \hspace{2cm}
+\nu\|u\|_{L^\infty_{s/2}(0,T; b^{2-s}_{q,\infty}(\Om))} +\|u_t, \na
p\|_{L^\infty_{s/2}(0,T; b^{-s}_{q,\infty}(\Om))}\leq M\nu^{1+s/2}
\end{array}
\end{equation}
for all $\theta\in (s/2,1)$. The solution $(u,\na p)$ is unique in
the class of functions satisfying \eq{3.6} with sufficiently small
$M=M(q,n,\Om,s)>0$ on the right-hand side.

}
\end{theorem}
{\bf Proof:} The proof is based on linearization and a fixed point
argument.

Consider the linear system \eq{2.21} with arbitrarily fixed $f\in
L^\infty_{s/2}(0,T;b^{-s}_{q,\infty}(\Om))$ and $u_0\in
B^0_{q,\infty,0,\si}(\Om)$. Then, by Theorem \ref{T2.7} the system
\eq{2.21} has a unique solution $\{u,\na p\}$ such that
$$\begin{array}{l}
u\in L^\infty(0,T; B^{0}_{q,\infty}(\Om))\cap L^\infty_{s/2}(0,T;
b^{2-s}_{q,\infty}(\Om)) \cap L^\infty_{1-\theta+s/2}(0,T;
B^{2\theta-s}_{q,1}(\Om)), \theta\in (\frac{s}{2},1),\ek
  u_t, \na p\in  L^\infty_{s/2}(0,T; b^{-s}_{q,\infty}(\Om))
\end{array}$$
 and
\begin{equation}
\label{E3.2}
\begin{array}{l}
 \nu^{s/2}\|u\|_{L^\infty(0,T; B^{0}_{q,\infty}(\Om))}
 +\nu^{\theta}\|u\|_{L^\infty(0,T; B^{2\theta-s}_{q,1}(\Om))}
  +\nu \|u\|_{L^\infty_{s/2}(0,T; b^{2-s}_{q,\infty}(\Om))}
  \ek
\hspace{0.5cm} + \|u_t, \na p\|_{L^\infty_{s/2}(0,T; b^{-s}_{q,\infty}(\Om))}
    \leq \tilde{C}(\|f\|_{L^\infty_{s/2}(0,T;b^{-s}_{q,\infty}(\Om))}
    +\nu^{s/2}\|u_0\|_{B^{0}_{q,\infty}(\Om)})
\end{array}
\end{equation}
with constant $\tilde{C}>0$ depending on $q,n,\Om$ and $s$. Now,
given $u_0\in B^{0}_{q,\infty,0,\si}(\Om)$, define the mapping
$\Phi$ from $Y:=L^\infty_{s/2}(0,T;b^{-s}_{q,\infty}(\Om))$ to
itself by
$$\Phi f:= -(u_f\cdot\na )u_f + a (\nu\Da u_f-\na p_f),$$
where $(u_f,\na p_f)$ is the unique solution to \eq{2.21}
corresponding to $u_0$ and $f$. If we show that $\Phi$ has a fixed
point $\tilde{f}\in Y$, then
 $(u_{\tilde{f}},\na p_{\tilde{f}})$ is obviously a solution to the system \eq{3.1}.

For $f_1,f_2\in Y$ and almost all $t\in (0,T)$ we get,  using
\eq{3.2} and Lemma \ref{L3.0} in view of $H^{-s}_{q}(\Om)
\hookrightarrow b^{-s}_{q,\infty}(\Om)$, that
 $$\begin{array}{rcl}
 \|(u_{f_1}\cdot\na)u_{f_2}(t)\|_{b^{-s}_{q,\infty}(\Om)} &
    \leq &c \|(u_{f_1}\cdot\na)u_{f_2}(t)\|_{H^{-s}_{q}(\Om)}\ek
    &\leq &c\|u_{f_1}(t)\|_{H^{2(1-\theta)}_{q}(\Om)} \|\na u_{f_2}(t)\|_{H^{2\theta-s-1}_{q}(\Om)}\ek
    &\leq  &c \|u_{f_1}(t)\|_{B^{2(1-\theta)}_{q,1}(\Om)} \|u_{f_2}(t)\|_{B^{2\theta-s}_{q,1}(\Om)}\ek
    &\leq  &c \nu^{-1-s/2}t^{-1+s/2}(\|{f_1}\|_{L^\infty_{s/2}(0,T;b^{-s}_{q,\infty}(\Om))}
             +\nu^{s/2}\|u_0\|_{B^{0}_{q,\infty}(\Om)})\cdot\ek
&&
\hspace{2.5cm}\cdot(\|{f_2}\|_{L^\infty_{s/2}(0,T;b^{-s}_{q,\infty}(\Om))}
             +\nu^{s/2}\|u_0\|_{B^{0}_{q,\infty}(\Om)}),
 \end{array}$$
 yielding
\begin{equation}
\label{E3.3}
\begin{array}{l}
\|(u_{f_1}\cdot\na)u_{f_2}\|_Y \leq
c\nu^{-1-s/2}(\|{f_1}\|_{L^\infty_{s/2}(0,T;b^{-s}_{q,\infty}(\Om))}
   +\nu^{s/2}\|u_0\|_{B^{0}_{q,\infty}(\Om)})\cdot\ek
\hspace{3cm}\cdot
(\|{f_2}\|_{L^\infty_{s/2}(0,T;b^{-s}_{q,\infty}(\Om))}
+\nu^{s/2}\|u_0\|_{B^{0}_{q,\infty}(\Om)})
\end{array}
\end{equation}
with $c=c(q,n,s,\Om)$. Hence, using \eq{3.2} and \eq{3.3}, we have
$$\begin{array}{rl}
\|\Phi(f)\|_Y&\leq \|(u_f\cdot\na)u_f\|_Y+\|a(\nu\Da u_f-\na
p_f)\|_Y\ek
 &\leq
C_0(\|f\|_{L^\infty_{s/2}(0,T;b^{-s}_{q,\infty})}
+\nu^{s/2}\|u_0\|_{B^{0}_{q,\infty}})
 (\|a\|_{L^\infty(0,T;{\cal M}(b^{-s}_{q,\infty}))}\ek
&\hspace{3cm}
 +\nu^{-1-s/2}\|f\|_{L^\infty_{s/2}(0,T;b^{-s}_{q,\infty})}
 +\nu^{-1}\|u_0\|_{B^{0}_{q,\infty}}),
 \end{array}$$
where $C_0=C_0(q,n,s,\Om)$.
Therefore, if $f\in B_{K}$,
where $B_K$ is the ball of $Y$ centered at $0$ with radius $K>0$,
then
\begin{equation}
\label{E3.9}
\|\Phi(f)\|_Y \leq C_0(K+\nu^{s/2}\|u_0\|_{B^{0}_{q,\infty}})
 (\|a\|_{L^\infty(0,T;{\cal M}(b^{-s}_{q,\infty}))}
 +\nu^{-1-s/2}K+\nu^{s/2}\|u_0\|_{B^{0}_{q,\infty}}).
\end{equation}

On the other hand, if $f_1,f_2\in B_{K}$, then
\begin{equation}
\label{E3.5}
\begin{array}{l}
\Phi(f_1)-\Phi(f_2)
=[-((u_{f_1}-u_{f_2})\cdot\na)u_{f_1}-(u_{f_2}\cdot\na)(u_{f_1}-u_{f_2})]
\ek\hspace{4cm} +[a(\nu\Da (u_{f_1}-u_{f_2})-\na (p_{f_1}-p_{f_2}))]
   =:\text{(I)+(II)}.
\end{array}
\end{equation}
Note that $u_f-u_g$ is a solution to \eq{2.21} with zero initial
value and right-hand side ${f_1}-{f_2}$. Hence we get by \eq{3.3}
that
$$\|\text{(I)}\|_Y\leq C_0\nu^{-1-s/2}\|f_1-f_2\|_Y (2K+\nu^{s/2}\|u_0\|_{B^{0}_{q,\infty}})$$
and, by \eq{3.2},
$$\|\text{(II)}\|_Y\leq C_0\|a\|_{L^\infty(0,T;{\cal M}(b^{-s}_{q,\infty}))}\|f_1-f_2\|_Y,$$
where we may regard the constant  $C_0$ exactly the same as in
\eq{3.9} in view of the structure of $\Phi(f_1)-\Phi(f_2)$, see
\eq{3.5}. Finally, we have
\begin{equation}
\label{E3.11} \|\Phi(f_1)-\Phi(f_2)\|_Y\leq C_0
    (2\nu^{-1-s/2}K+\|a\|_{L^\infty(0,T;{\cal M}(b^{-s}_{q,\infty}))}
    +\nu^{-1}\|u_0\|_{B^{0}_{q,\infty}}
 )\|f_1-f_2\|_Y.
\end{equation}

Now, in view of \eq{3.9} and \eq{3.11}, consider the inequality
\begin{equation}
\label{E3.12}
\left\{\begin{array}{l}
              C_0(K+\nu^s\|u_0\|_{B^{0}_{q,\infty}})(\nu^{-1-s/2}K
              +\|a\|_{L^\infty(0,T;{\cal M}(b^{-s}_{q,\infty}))}
                 +\nu^{-1}\|u_0\|_{B^{0}_{q,\infty}})<K,\ek
              C_0 (2\nu^{-1-s/2}K+\|a\|_{L^\infty(0,T;{\cal M}(b^{-s}_{q,\infty}))}
                +\nu^{-1}\|u_0\|_{B^{0}_{q,\infty}})<1.
         \end{array}
  \right.
\end{equation}
 By elementary calculations, it follows that,
if
\begin{equation}
\label{E3.4}
\|a\|_{{\cal M}(b^{-s}_{q,\infty})}<\frac{3}{8C_0},
 \quad\nu^{-1}\|u_0\|_{B^{0}_{q,\infty}}<\frac{1}{8C_0},
\end{equation}
then for any
\begin{equation}
\label{E3.10}
K\in (k_1,k_2)
\end{equation}
with
$$k_{1,2}=\frac{\frac{1}{2}-\nu^{-1}\|u_0\|C_0\mp
       \sqrt{(\frac{1}{2}-\nu^{-1}\|u_0\|C_0)^2-2C_0\nu^{-1}
\|u_0\|_{B^{0}_{q,\infty}}}}{2\nu^{-1-s/2}C_0}$$
  the inequality
\eq{3.12} holds true. In other words, if \eq{3.4} and \eq{3.10} are
satisfied, then $\Phi(B_{K})\subset B_{K}$ and $\Phi: B_{K}\mapsto
B_{K}$ is a contraction mapping. Thus, by the Banach fixed point
theorem $\Phi$ has a fixed point $\tilde{f}$ in $B_K$, which is
unique in $B_{K}$, and $u=u_{\tilde{f}}$ is a solution to \eq{3.1}.

Note that $K<\frac{\nu^{1+s/2}}{2C_0}$.
It follows from \eq{3.2}, \eq{3.4} and \eq{3.10} that the solution $u$
satisfies \eq{3.6} with $M\equiv \frac{5\tilde{C}}{8C_0}$.

Let $K'\equiv \frac{K}{3}$ and the norm of solution $(u, \na p)$ in
\eq{3.6} is bounded by $K'$. Then $\|u_t-\nu\Da u+\na
p\|_{L^\infty_\ga(0,T; b^{-s}_{q,\infty})}\leq K$. Hence, in view of
the uniqueness of the linear Stokes problem, $(u,\na p)$ must be the
only solution satisfying the inequality \eq{3.6} with $K'$ on the
right-hand side.

The proof of the theorem is complete. \qed
%
%

\section{Transport equation}
In this section, we shall prove an existence result for the
transport equation with piecewise constant initial values.

We know that the characteristic function $\chi(\Om')$ of any
Lipschitz subset $\Om'$ of $\Om$ can be a pointwise multiplier of
Besov spaces $B^s_{q,r}(\Om)$ and $b^s_{q,\infty}(\Om)$,
respectively, with $q\in (1,\infty), r\in [1,\infty)$, $s\in
(-1+1/q,1/q)$, cf. \cite{Tri02}, Proposition 5.1, 5.3; cf. also
\cite{Ya10}, Theorem 1.41. Here, the point is that for $q,r,s$
satisfying the above conditions the set $C^\infty_0(\Om')$ is dense
in $B^s_{q,r}(\Om')$ and $b^s_{q,\infty}(\Om')$, thus extension of
$f\in B^s_{q,r}(\Om')$ or $f\in b^s_{q,\infty}(\Om')$ by $0$ in
$\Om\setminus\Om'$ defines a linear continuous extension with norm
$1$ from  $B^s_{q,r}(\Om')$ to $B^s_{q,r}(\Om)$ and
$b^s_{q,\infty}(\Om')$ to $b^s_{q,\infty}(\Om)$, respectively.


\begin{lemma}
\label{L3.1} {\rm Let $\Om$ be a domain of $\R^n$, $n\in\N$, and let
$1<q<\infty$ and $s\in (-1+1/q,1/q)$. Then, for any Lipschitz
subdomain $\Om'$ of $\Om$
 $$\|\chi_{\Om'}u\|_{Y}\leq c\|u\|_{Y},\quad\forall u\in Y,$$
with $c>0$ independent of $\Om'$, where $Y=B^{s}_{q,r}(\Om)$, $1\leq
r<\infty$, or $Y=b^s_{q,\infty}(\Om)$.
 }
\end{lemma}
{\bf Proof:} First, suppose that $s\in (0,1/q)$.
Let $\tilde{L}^q(\Om)$ and $\tilde{H}^1_q(\Om)$ be, respectively,
the (closed) subspace of $L^q(\Om)$ and $H^1_q(\Om)$ of functions with
support in $\Om'$.
 Define the operator $E_0: L^q(\Om')\mapsto \tilde{L}^q(\Om)$ by $E_0f:= \tilde{f}$,
 where $\tilde{f}$ is the extension of $f$ by zero on $\Om\setminus\Om'$.
Then,
$$E_0\in {\cal L}(L^q(\Om'), \tilde{L}^q(\Om))\cap {\cal L}(H^1_{q,0}(\Om'), \tilde{H}^1_q(\Om))$$
and
$$\|E_0\|_{{\cal L}(L^q(\Om'), \tilde{L}^q(\Om))}=1, \quad
 \|E_0\|_{{\cal L}(H^1_{q,0}(\Om'), \tilde{H}^1_q(\Om))}=1.$$
Hence, by real interpolation $(\cdot,\cdot)_{s,r}$, $1\leq
r<\infty$, and $(\cdot,\cdot)^0_{s,\infty}$,
 we have
$$E_0\in {\cal L}\big((L^q(\Om'), H^1_{q,0}(\Om'))_{s,r},
(\tilde{L}^q(\Om), \tilde{H}^1_q(\Om))_{s,r}\big), 1\leq r<\infty,$$
and
$$E_0\in {\cal L}\big((L^q(\Om'), H^1_{q,0}(\Om'))^0_{s,\infty},
(\tilde{L}^q(\Om), \tilde{H}^1_q(\Om))^0_{s,\infty}\big),
\quad\|E_0\|\leq 1.$$ Note that, due to $s\in (0,1/q)$, we have
$$(L^q(\Om'), H^1_{q,0}(\Om'))_{s,r}=B^s_{q,r}(\Om'),\quad
(L^q(\Om'), H^1_{q,0}(\Om'))^0_{s,\infty}=b^s_{q,\infty}(\Om'),$$
see \cite{Am00}, Theorem 2.2. Moreover,  by interpolation, we have
continuous embedding
$$\begin{array}{l}
(\tilde{L}^q(\Om), \tilde{H}^1_q(\Om))_{s,r}\subset (L^q(\Om), H^1_q(\Om))_{s,r}=B^s_{q,r}(\Om),\ek
(\tilde{L}^q(\Om), \tilde{H}^1_q(\Om))^0_{s,\infty}\subset
(L^q(\Om), H^1_q(\Om))^0_{s,\infty}=b^s_{q,\infty}(\Om)
\end{array}$$
with embedding constants not greater than $1$.
Therefore,
\begin{equation}
\label{E2.3}
E_0\in {\cal L}\big(B^s_{q,r}(\Om'),B^s_{q,r}(\Om)\big),\quad
\|E_0\|_{{\cal L}\big(B^s_{q,r}(\Om'),B^s_{q,r}(\Om)\big)}\leq 1
\end{equation}
and
\begin{equation}
\label{E2.4}
E_0\in {\cal L}\big(b^s_{q,\infty}(\Om'),b^s_{q,\infty}(\Om)\big),\quad
\|E_0\|_{{\cal L}\big(b^s_{q,\infty}(\Om'),b^s_{q,\infty}(\Om)\big)}\leq 1.
\end{equation}

On the other hand, for $r_{\Om'}$ being the restriction operator
onto $\Om'$, we have
$$\begin{array}{l}
\|r_{\Om'}f\|_{L^{q}(\Om')}\leq \|f\|_{L^{q}(\Om)},\quad\forall f\in L^{q}(\Om), \ek
\|r_{\Om'}f\|_{H^1_{q}(\Om')}\leq \|f\|_{H^1_{q}(\Om)}, \quad\forall f\in H^1_{q}(\Om).
\end{array}$$
Hence, by real interpolation we have
\begin{equation}
\label{E2.5}
\begin{array}{l}
\|r_{\Om'}f\|_{B^s_{q,r}(\Om')}\leq \|f\|_{B^s_{q,r}(\Om)},\quad\forall f\in B^s_{q,r}(\Om), \ek
\|r_{\Om'}f\|_{b^s_{q,\infty}(\Om')}\leq \|f\|_{b^s_{q,\infty}(\Om)},
\quad\forall f\in b^s_{q,\infty}(\Om).
\end{array}
\end{equation}

Combining \eq{2.3}-\eq{2.5}, we get that
\begin{equation}
\label{E2.6}
\begin{array}{l}
\|E_0r_{\Om'}f\|_{B^s_{q,r}(\Om)}\leq \|f\|_{B^s_{q,r}(\Om)},
\quad\forall f\in B^s_{q,r}(\Om),\ek
\|E_0r_{\Om'}f\|_{b^s_{q,\infty}(\Om)}\leq \|f\|_{b^s_{q,\infty}(\Om)},
 \quad\forall f\in b^s_{q,\infty}(\Om).
\end{array}
\end{equation}
In \eq{2.3}-\eq{2.6}, $B^s_{q,r}(\Om)$, $r\in [1,\infty)$,
 is endowed with the norm of the
real interpolation space $(L^q(\Om), H^1_q(\Om))_{s,r}$ and
$b^s_{q,\infty}(\Om)$ is endowed with the
 norm
 of $B^s_{q,\infty}(\Om)$; both norms are equivalent to
 the typical Besov norms
 $$\|u\|_{B^s_{q,r}(\Om)}=\inf_{u=r_\Om \tilde{u}, \tilde{u}\in B^s_{q,r}(\R^n)}
 \|\tilde{u}\|_{B^s_{q,r}(\R^n)}\,\,\text{and }
\|u\|_{b^s_{q,\infty}(\Om)}=\inf_{u=r_\Om \tilde{u}, \tilde{u}\in
b^s_{q,\infty}(\R^n)}
 \|\tilde{u}\|_{B^s_{q,\infty}(\R^n)},$$
 respectively.
Thus, in view of $\chi_{\Om'}f=E_0r_{\Om'}f$,
we get the assertion of the lemma for $s\in (0,1/q)$.

The assertion of the lemma for the case $s\in (-1+1/q,0)$ follows by duality
argument using the result for $s\in (0,1/q)$.

Finally, the case $s=0$ directly follows by interpolation. \qed
\par\bigskip
Based on Lemma \ref{L3.1}, we can prove the following statement.
\begin{proposition}
\label{P4.2} {\rm Let $\Om\subset \R^n,n\geq 3,$ be whole or half
space, or bounded domain with boundary of $C^2$-class, and let $u\in
L^1(0,T; W^{1,\infty}(\Om))$, $\div u=0$ and $u|_{\pa\Om}=0$. Let
$\Om_1$ be a Lipschitz subdomain of $\Om$ and $\Om_2=\Om\setminus
\bar\Om_1$ and
$\rho_0(x)=\rho_{01}\chi_{\Om_1}(x)+\rho_{02}\chi_{\Om_2}(x)$,
$x\in\Om$,
 $ 0<\rho_{01}<\rho_{02}$.
Then, the transport equation
\begin{equation}
\label{E3.25}
\rho_t+u\cdot\na\rho=0\quad\text{in }(0,T)\ti\Om,\quad \rho(0)=\rho_0\quad\text{in }\Om,
\end{equation}
has a unique solution $\rho$ such that
for all $q\in (1,\infty)$, $s\in (-1+1/q,1/q)$
\begin{equation}
\label{E3.26} \rho \in L^\infty(0,T;L^\infty(\Om))\cap
L^\infty(0,T;{\cal M}(Y))
\end{equation}
and
\begin{equation}
\label{E3.27}
\begin{array}{l}
\|\rho(t)\|_{L^\infty(\Om)}= \|\rho_0\|_{L^\infty(\Om)}, \,\,\forall
t\in (0,T),\ek
 \|\rho\|_{L^\infty(0,T;{\cal
M}(Y))}\leq c\rho_{02},\quad \|a\|_{L^\infty(0,T;{\cal M}(Y))}\leq
c(\frac{\bar\rho-\rho_{01}}{\rho_{01}}
+\frac{\rho_{02}-\bar\rho}{\rho_{02}}),
\end{array}
\end{equation}
where $Y=B^{s}_{q,r}(\Om)$, $1\leq r<\infty$, or
$Y=b^s_{q,\infty}(\Om)$, $a=\frac{\bar\rho}{\rho}-1$ and
$c=c(q,s,\Om)$.

}
\end{proposition}
{\bf Proof:} Under the assumption of the proposition,
 unique existence of solution in $L^\infty(0,T;L^\infty(\Om))$
is already proved in \cite{DiLi89}, Theorem II.3, more precisely,
the solution $\rho$ to \eq{3.25} is expressed by
$\rho(t,x)=\rho_0(X(t,\cdot)^{-1}(x))$ and
$\|\rho(t)\|_{L^\infty(\Om)}=\|\rho_0\|_{L^\infty(\Om)}$, $\forall
t\in (0,T)$, where $\{X(t,y)\}_{t\geq 0}$ stands for the semiflow
transported by the vector field $u$, that is,
$$X(t,y)=y+\int_0^t u(\tau, X(\tau,y))\,d\tau, \quad t\in (0,T), y\in \Om.$$
Note that $X(t,\cdot)$  is  a $C^1$-diffeomorphism
over $\Om$ for each $t>0$.


Since $\rho_0(x)=\rho_{01}\chi_{\Om_1}(x)+\rho_{02}\chi_{\Om_2}(x)$,
we have
\begin{equation}
\label{E3.29}
\begin{array}{l}
\rho(t,x)=\rho_{01}\chi_{\Om_{1}(t)}(x)+\rho_{02}\chi_{\Om_{2}(t)}(x),\ek
a(t,x)=\frac{\bar\rho-\rho_{01}}{\rho_{01}}\chi_{\Om_1(t)}(x)
+\frac{\rho_{02}-\bar\rho}{\rho_{02}}\chi_{\Om_2(t)}(x), t\in (0,T),
x\in\Om,
\end{array}
\end{equation}
where
\begin{equation}
\label{E3.7} \Om_{i}(t)=\{X(t,y): y\in \Om_{i}\}, i=1,2.
\end{equation}
Here we used that $\Om_1(t)\cap\Om_2(t)=\emptyset$ for all $t\in
(0,T)$ due to Lipschitz conditions on the vector field $u$.
Therefore we get \eq{3.27} by Lemma \ref{L3.1}.

\qed
\begin{remark} {\rm
\label{R4.3} In Lemma \ref{L3.1} the boundary of $\Om'$ is allowed
to be {\it so-called} a $d$-set with some $d\in (0,n)$ and
\begin{equation}
\label{E4.31} s\in \big((n-d)/q', (n-d)/q\big).
\end{equation}
We recall the definition of a $d$-set for
$0<d\leq n$; a closed non-empty set $\Gamma$ of $\R^n$ is called a
$d$-set if
$$\exists c_1, c_2>0: \forall x\in\Gamma,\forall r\in (0,1],
c_1r^d\leq {\cal H}^d (B(x,r)\cap \Gamma)\leq c_2 r^d, $$
 where ${\cal H}^d$ denotes the $d$-dimensional Hausdorff measure on $\R^n$ and
$B(x,r)$ stands for the open ball centered at $x$ with radius $r$
(cf. e.g. \cite{Ca00}). In fact, if $\Om'\subset \R^n$ is a $d$-set
with some $d\in (0,n)$, then, by \cite{Ca00}, Corollary 2.7 it
follows that $C^\infty_0(\Om')$ is dense in $B^s_{q,r}(\Om')$,
$1<q,r<\infty$ for all $s\in (0, (n-d)/q)$ and hence, by duality,
for all $s$ satisfying \eq{4.31}. Consequently, by denseness of
$B^s_{q,r}(\Om')$ in $b^s_{q,\infty}(\Om')$,
 we get that
$C^\infty_0(\Om')$ is dense in $b^s_{q,\infty}(\Om')$.

Therefore, in Proposition \ref{P4.2}, the initial interface of two
fluids $\pa\Om_1\cap \pa\Om_2$ is also allowed to be a $d$-set. Note
that, if
  $X(t,\cdot)$, $t\in (0,T)$,
 is diffeomorphism over $\Om$, the boundary $\pa\Om_i(t)$, $i=1,2$, of $\Om_i(t)$, $i=1,2$,
 given by \eq{3.7} remains as $d$-sets for all time $t\in [0,T)$  due to the property
 $${\cal H}^d(f({\cal A}))\leq (Lip (f))^d {\cal H}^d({\cal A}), \quad {\cal A}\subset\R^n,$$
for Lipschitz function $f$, where $Lip (f)$ is the Lipschitz
constant of $f$, see e.g. \cite{EvaGa92}, Section 2.4, Theorem 1.

}
\end{remark}

\section{Proof of Theorem \ref{T1.2}}
The procedure to prove Theorem \ref{T1.2} is twofold, i.e.,
existence part and uniqueness part.

\subsection{Proof of existence}
Let $\Om$ be a bounded domain with $C^2$-boundary of $\R^n$, $n\geq
3$, and let $\Om_1$ be a subdomain of $\Om$ with Lipschitz boundary
and $\Om_2=\Om\setminus \bar\Om_1$. Suppose that
\begin{equation}
\label{E4.3}
\rho_0(x)=\rho_{01}\chi_{\Om_1}(x)+\rho_{02}\chi_{\Om_2}(x),
\;x\in\Om,\; 0<\rho_{01}<\rho_{02}, \quad u_0\in
B^0_{q,\infty,0,\si}(\Om), q\geq n.
\end{equation}

Let $\eta_m\in C^\infty(\R^n), m\in\N,$ be mollifiers
such that
$$\eta_m(x)\geq 0,\,\,\eta(x)=\eta(-x),\,\,
\supp \eta_m\subset \big\{x: |x| <\frac{1}{m}\big\},\,\,
 \int_{\R^n} \eta_m(x)\,dx=1.$$

For $m=1,2,\ldots$, let us construct an iterate scheme for \eq{1.1}
as
\begin{equation}
\label{E4.4}
\left\{
\begin{array}{l}
     \rho_{mt}+(u^{(m-1)}\cdot\na)\rho_m=0,\,\,\rho_{m}(0,x)=\rho_0(x),
      \quad\text{in }(0,T)\ti\Om,\ek
 u_{mt}-\nu\Da u_m  + (u_m \cdot \na) u_m+\na p_m  = a_m(\nu\Da u_m-\na p_m),
 \quad \text{in }(0,T)\ti\Om,\ek
      \div u_m = 0,\quad \text{in } (0,T)\ti\Om,\ek
     u_m = 0, \quad \text{on }(0,T)\ti\pa\Om,\ek
  u_m(0,x) = u_0, \quad\text{in }\Om,
\end{array}
\right.
\end{equation}
where $\nu:=\frac{\mu}{\bar\rho}$ with $\bar\rho\in
(\rho_{01},\rho_{02})$ fixed, $0<T\leq \infty$, $u^{(0)}\equiv 0$
  and $a_m(t,x):=\frac{\bar\rho}{\rho_m(t,x)}-1$,
 $u^{(m)}=\eta_m\star\bar{u}_m|_\Om$,
$$\bar{u}_m=\left\{
           \begin{array}{cl}
            u_m &   \text{for }x\in \Om\,\,\text{with dist}(x,\pa\Om)>\frac{1}{m}\ek
             0  &   \text{for else }x\in \R^n
          \end{array}
          \right.
$$
for $m\in\N$. Obviously, $u^{(m)}|_{\pa\Om}=0$, $m\in\N$.

\begin{remark}
\label{R4.1} {\rm Let $\Om$ be a bounded domain of $\R^n,n\in\N,$
with Lipschitz boundary. Denoting  by $\tilde{w}$ the  extension of
$w$ by zero in $\R^n\setminus \Om$, for $w\in b^s_{q,\infty}(\Om)$,
$q\in (1,\infty)$, $s\in (-1+1/q,1/q)$ one gets $\tilde{w}\in
b^s_{q,\infty}(\R^n)$ and $\|\tilde{w}\|_{b^s_{q,\infty}(\R^n)}\leq
c(q,s,\Om)\|w\|_{b^s_{q,\infty}(\Om)}$.

On the other hand, for $f\in  H^1_q(\R^n)$ (or $H^{-1}_q(\R^n)$) it
holds $\|\eta_m\star f-f\|_{H^1_q(\R^n)}\ra 0$ (or $\|\eta_m\star
f-f\|_{H^{-1}_q(\R^n)}\ra 0$ ) as $m\ra\infty$. Hence, by the
Banach-Steinhaus theorem one gets uniform boundedness  with respect
to $m\in\N$ of the operator norms of convolution operators
$\eta_m\star\cdot$ in ${\cal L}(H^1_q(\R^n))\cap {\cal
L}(H^{-1}_q(\R^n))$ and, consequently, in ${\cal
L}(B^{1-2\theta}_{q,\infty}(\R^n))\cap {\cal
L}(b^{1-2\theta}_{q,\infty}(\R^n))$, $\theta\in (0,1)$, due to
$(H^1_q(\R^n),H^{-1}_q(\R^n))_{\theta,\infty}=B^{1-2\theta}_{q,\infty}(\R^n)$
and
$(H^1_q(\R^n),H^{-1}_q(\R^n))^0_{\theta,\infty}=b^{1-2\theta}_{q,\infty}(\R^n)$.

 Therefore it follows  that
 \begin{equation}
\label{E4.1}
\begin{array}{l}
(w\mapsto(\eta_m\star\tilde{w})|_\Om)\in {\cal L}(b^s_{q,\infty}(\Om),b^s_{q,\infty}(\R^n)),
\forall s\in (-1+1/q,1/q),\ek
\|(\eta_m\star \tilde{w})|_\Om\|_{b^s_{q,\infty}(\Om)}\leq
 C\|w\|_{b^s_{q,\infty}(\Om)},\,\,
\forall w\in b^s_{q,\infty}(\Om),
\end{array}
\end{equation}
with $C=C(q,s,\Om)>0$ independent of $m\in\N$.

Moreover,
it follows that
\begin{equation}
\label{E4.5}
(\eta_m\star \tilde{w})|_\Om\ra w\,\,\text{in }b^s_{q,\infty}(\Om),\,\,
 \forall w\in b^s_{q,\infty}(\Om), -1+1/q< s<1/q,\quad(\text{as }m\ra\infty).
\end{equation}
In fact, if $w\in H^1_{q,0}(\Om)$, then $(\eta_m\star
\tilde{w})|_\Om$ tends to $w$ in $H^1_{q,0}(\Om)$ and
$H^{-1}_q(\Om)$, respectively, hence in $b^s_{q,\infty}(\Om)$, $-1<
s<1$, by continuous interpolation. Thus we get \eq{4.5}, in view of
\eq{4.1} and denseness of $H^1_{q,0}(\Om)$
  in $b^s_{q,\infty}(\Om)$, $-1< s<1/q$.

Furthermore, it follows that, if $u_m\in L^\infty_{s/2}(0,T;
b^{2-s}_{q,\infty}(\Om))$, $q\geq n$, for some $s\in (0,2)$,
 then
$u^{(m)}(t)\in C^\infty(\R^n)$, $\supp {u^{(m)}(t)}\subset\bar\Om$
for almost all $t\in (0,T),$ and, in particular,
\begin{equation}
\label{}
u^{(m)}\in L^1_{loc}(0,T; W^{1,\infty}(\Om)).
\end{equation}
%
 }
\end{remark}
%
%
\par\bigskip
We have the following lemma.

\begin{lemma}
\label{L4.3} {\rm Let $\Om$ be a bounded domain of $C^2$-class and
suppose that \eq{4.3} holds. Then
 for any $s\in (0,1-\frac{1}{q})$
 there are some constants $\da_i=\da_i(q,n,s,\Om)>0, i=1,2,$ and
$M=M(q,n,s,\Om)>0$ independent of $m\in\N$ with the following
property: If $$\frac{\rho_{02}-\rho_{01}}{\rho_{01}}<\da_1,\quad
\|u_0\|_{B^0_{q,\infty}(\Om)}<\da_2\nu,$$
 the system \eq{4.4} has a solution
$\{(\rho_m,u_m, p_m):m\in\N\}$ satisfying
\begin{equation}
\label{E4.24}
\begin{array}{l}
\rho_m \in L^\infty(0,T;L^\infty(\Om))\cap L^\infty(0,T;{\cal
M}(b^{\tilde{s}}_{\tilde{q},\infty}(\Om))),\ek
 \|\rho_m\|_{L^\infty(0,T;L^\infty(\Om))}=\|\rho_0\|_{L^\infty(\Om)},\ek
  \|\rho_m\|_{L^\infty(0,T;{\cal
M}(Y))} \leq c(\tilde{q},\tilde{s},\Om)\rho_{02},\ek
 \|a_m\|_{L^\infty(0,T;{\cal M}(Y))} \leq
c(\tilde{q},\tilde{s},\Om)(\frac{\bar\rho-\rho_{01}}{\rho_{01}}
+\frac{\rho_{02}-\bar\rho}{\rho_{02}}),
\end{array}
\end{equation}
where $Y=B^{\tilde{s}}_{\tilde{q},r}(\Om)$, $1\leq r<\infty$, or
$Y=b^{\tilde{s}}_{\tilde{q},\infty}(\Om)$
 for any $\tilde{q}\in (1,\infty)$, $\tilde{s}\in (-1+1/\tilde{q},1/\tilde{q})$, and
\begin{equation}
\label{E4.25}
\begin{array}{l}
u_m\in L^\infty(0,T;B^{0}_{q,\infty}(\Om))\cap
  L^\infty_{s/2}(0,T; b^{2-s}_{q,\infty}(\Om))
  \cap L^\infty_{1-\theta+s/2}(0,T; B^{2\theta-s}_{q,1}(\Om)),\ek
u_{mt}, (u_m\cdot\na)u_m, \na p_m\in L^\infty_{s/2}(0,T; b^{-s}_{q,\infty}(\Om))
\end{array}
\end{equation}
with estimate
\begin{equation}
\label{E4.26}
\begin{array}{l}
\nu^{s/2}\|u_m\|_{L^\infty(0,T; B^{0}_{q,\infty}(\Om))}
+\nu\|u_m\|_{L^\infty_{s/2}(0,T;
b^{2-s}_{q,\infty}(\Om))}\ek\hspace{5cm} +\|u_{mt}, \na
p_m\|_{L^\infty_{s/2}(0,T; b^{-s}_{q,\infty}(\Om))}\leq
M\nu^{1+s/2},\ek \|u_m\|_{L^\infty_{1-\theta+s/2}(0,T;
B^{2\theta-s}_{q,1}(\Om))}\leq M\nu^{1+s/2-\theta},\quad\forall
\theta\in (s/2,1).
\end{array}
\end{equation}

 }
\end{lemma}
{\bf Proof:} This lemma follows directly by Theorem \ref{T2.4} and
Proposition \ref{P4.2}, \qed

\begin{remark}
\label{R4.4} {\rm Let $v\in L^\infty_{1-\theta+s/2}(0,T;
B^{2\theta-s}_{q,1}(\Om))$ for all $\theta\in (s/2,1)$. Then,
choosing $\theta_1\in (s/2,1)$ for any $q_1\in (1,\frac{q(n+2)}{n})$
as
\begin{equation}
\label{E4.19}
\frac{s}{2}+\frac{1}{2}\big(\frac{n}{q}-\frac{n}{q_1}\big) <
\theta_1<\frac{s}{2}+\frac{n}{(n+2)q},
\end{equation}
we have
$$\frac{q(n+2)}{n}(\theta_1-\frac{s}{2})<1, \quad B^{2\theta_1-s}_{q,1}(\Om)
\hookrightarrow L^{q_1}(\Om).$$ Hence, in view of
$\|v(t)\|_{B^{2\theta_1-s}_{q,1}(\Om)}\leq ct^{-2\theta_1+s}$,
$\forall t\in (0,T)$, we have
\begin{equation}
\label{E4.12} v\in L^{q(n+2)/n}(0,T;B^{2\theta_1-s}_{q,1}(\Om))
\hookrightarrow
 L^{q_1}(Q_T)
 \end{equation}
and
\begin{equation}
\label{E4.7} \|v\|_{L^{q(n+2)/n}(0,T;B^{2\theta_1-s}_{q,1}(\Om))}
 \leq c\|v\|_{L^\infty_{1-\theta_1+s/2}((0,T), B^{2\theta_1-s}_{q,1}(\Om))}.
\end{equation}
Similarly, one gets for any $q_2\in (1,\frac{q(n+2)}{n+q})$ that
\begin{equation}
\label{E4.12n} \na v\in L^{q_2}(Q_T),\; \|\na v\|_{L^{q_2}(Q_T)}\leq
    c\|v\|_{L^\infty_{1-\theta_2+s/2}((0,T), B^{2\theta_2-s}_{q,1}(\Om))},
   \end{equation}
where
$$\frac{s+1}{2}+\frac{1}{2}\big(\frac{n}{q}-\frac{n}{q_2}\big)
< \theta_2<\frac{s}{2}+\frac{n+q}{(n+2)q}.$$
 Note that $q\geq n>2n/(n+1)$ and hence one can choose $q_2$ as
 $$ \Big(\frac{q(n+2)}{n}\Big)'\equiv \frac{q(n+2)}{q(n+2)-n}<q_2<\frac{q(n+2)}{n+q}.$$
Therefore, in view of \eq{4.12} and \eq{4.12n}, one has

\begin{equation}
\label{E4.11} \na v\in (L^{q(n+2)/n-\da}(Q_T))',\quad
 \|\na v\|_{(L^{q(n+2)/n-\da}(Q_T))'}\leq
 c\|v\|_{L^\infty_{1-\theta_2+s/2}((0,T), B^{2\theta_2-s}_{q,1}(\Om))}
\end{equation}
for sufficiently small $\da>0$.
 }
\end{remark}

\begin{remark}
\label{R4.5} {\rm Let $q\geq n$ and $\{(\rho_m,u_m, p_m):m\in\N\}$
be a solution to the system \eq{4.4} whose existence is guaranteed
by Lemma \ref{L4.3}. Then, since $u_m\in
L^\infty_{1-\theta+s/2}(0,T; B^{2\theta-s}_{q,1}(\Om))$, we get by
\eq{4.1} of Remark \ref{R4.1} that
$$u^{(m)}\in L^\infty_{1-\theta+s/2}(0,T; B^{2\theta-s}_{q,1}(\Om)),
 s\in (0,1-1/q), \theta\in (s/2,1), m\in\N.$$
Therefore, using $\div u^{(m-1)}=0$, $u^{(m-1)}|_{\pa\Om}=0$,
$m\in\N$, we get that if $s\in (0, 1-1/q), \theta\in
(s/2,s/2+1/(2q))$, then for $\vp\in C^\infty(\bar\Om)$ and almost
all $t\in (0,T)$
$$\begin{array}{rcl}
|\int_\Om (u^{(m-1)}(t)\cdot\na)\rho_m(t)\vp\,dx|
 &=&|\int_\Om \rho_m(t) u^{(m-1)}(t)\cdot \na\vp\,dx|\ek
&\leq& \|u^{(m-1)}(t)\|_{B^{2\theta-s}_{q,1}(\Om)}
\|\rho_m(t)\na\vp\|_{b^{-2\theta+s}_{q',\infty}(\Om)}\ek
&\leq& \|u^{(m-1)}(t)\|_{B^{2\theta-s}_{q,1}(\Om)}
\|\rho_m(t)\|_{{\cal M}(b^{-2\theta+s}_{q',\infty}(\Om))}
\|\vp\|_{b^{1-2\theta+s}_{q',\infty}(\Om)};
\end{array}$$
here note that
$(B^{2\theta-s}_{q,1}(\Om))'=B^{-2\theta+s}_{q',\infty}(\Om)$ thanks
to $0<2\theta-s<1/q$ (cf. \cite{Tri83}, Theorem 4.3.2/1 and (2.6.2)
of \cite{Am95}, \S 2.6). Thus, in view of
$\rho_{mt}=-(u^{(m-1)}\cdot\na)\rho_m$, $u^{(m-1)}\in
L^\infty_{1-\theta+s/2}(0,T;B^{2\theta-s}_{q,1}(\Om))$ and
$\rho_m\in L^\infty(0,T; {\cal
M}(b^{-2\theta+s}_{q',\infty}(\Om)))$, we have $
t^{\theta-s/2}\rho_{mt}\in L^\infty(0,T; H^{-1+2\theta-s}_q(\Om))$
and, in particular,
\begin{equation}
\label{E4.13} \rho_{mt}\in L^p(0,T;
H^{-1+2\theta-s}_q(\Om)),\,\,\forall p\in (1,\frac{2}{2\theta-s}),
\end{equation}
where $s\in (0, 1-1/q), \theta\in (s/2,s/2+1/(2q))$.
 }
\end{remark}
\par\bigskip
{\bf Proof of Theorem \ref{T1.2}: existence part}

Let $s\in (0, 1/q)$ and let
$\{(\rho_m, u_m, p_m):m\in\N\}$ be the solutions to the iterate system \eq{4.4},
 the existence of which is given by Lemma \ref{L4.3}.
 Then, by Lemma \ref{L4.3} $\{u_m\}$ is
 bounded in
$L^\infty(0,T; B^{0}_{q,\infty}(\Om))\cap L^\infty_{s/2}(0,T;
b^{2-s}_{q,\infty}(\Om))$,
 $\{u_{mt}\}$ is
bounded in $L^\infty_{s/2}(0,T; b^{-s}_{q,\infty}(\Om))$ and
$\{\rho_m\}$ is bounded in $L^\infty(Q_T)$. Hence, $\{u_m\}$ and
$\{\rho_m\}$ have some subsequences $\{u_{m_k}\}$ and
$\{\rho_{m_k}\}$, respectively, such that
\begin{equation}
\label{E4.8}
\begin{array}{ll}
u_{m_k} \rightharpoonup u &\quad \text{in }
           L^\infty(0,T; B^{0}_{q,\infty}(\Om)^{n})\,\,(\text{*-weakly as }k\ra\infty),\ek
u_{m_k} \rightharpoonup u &\quad \text{in }
           L^\infty_{s/2}(0,T; b^{2-s}_{q,\infty}(\Om)^{n})
           \,\,(\text{*-weakly as }k\ra\infty),\ek
u_{{m_k}t} \rightharpoonup u_{t} &\quad \text{in }
           L^\infty_{s/2}(0,T; b^{-s}_{q,\infty}(\Om)^n)\,\,(\text{*-weakly as }k\ra\infty),\ek
\rho_{m_k} \rightharpoonup \rho &\quad \text{in }
           L^\infty(0,T; L^\infty(\Om))\,\,(\text{*-weakly as }k\ra\infty)
\end{array}
\end{equation}
for some $u$, $\rho$. Moreover, it follows from \eq{4.12}, \eq{4.7}
that
\begin{equation}
\label{E4.15}
u_{m_k} \rightharpoonup u \quad \text{in }
           L^{q(n+2)/n}(0,T; H^{2\theta_1-s}_{q}(\Om)^{n})
           \,\,(\text{weakly as }k\ra\infty)
\end{equation}
for all $\theta_1\in (s/2,1)$ satisfying \eq{4.19}.

We shall show that $(\rho, u)$ is a solution to \eq{1.1} in the
sense of Definition \ref{D1.1}. Note that by \eq{4.25}, \eq{4.26}
the sequence $\{u_{m_kt}\}$ weakly converges in $L^{\al}(0,T;
H^{-s-\ve}_q(\Om))$ for some $\al>1$ and any $\ve>0$. Therefore, in
view of \eq{4.15} and compact embedding
$H^{2\theta_1-s}_{q}(\Om)\hra\hra L^{q(n+2)/n-\da}(\Om)$ for
sufficiently small $\da>0$, we get by compactness theorem
(\cite{Te77}, Ch.3, Theorem 2.1) that
\begin{equation}
\label{E4.20} u_{m_k}\ra u\quad\text{in } L^{q(n+2)/n-\da}(Q_T)
\end{equation}
 as $k\ra\infty$.
 Moreover,  we get from \eq{4.11} that $\{\na u_{m_k}\}$ is bounded in
$(L^{q(n+2)/n-\da}(Q_T))'$.

Rewriting the second equation of the system \eq{4.4}, we have
$$u_{m_kt}-\frac{\mu}{\rho_{m_k}}\Da u_{m_k}+(u_{m_k}\cdot\na)u_{m_k}
+\frac{\bar\rho}{\rho_{m_k}}\na p_{m_k}=0,$$
which is equivalent to
\begin{equation}
\label{E4.14}
\rho_{m_k}u_{m_kt}-\mu\Da u_{m_k}+\rho_{m_k}(u_{m_k}\cdot\na)u_{m_k}+\bar\rho\na p_{m_k}=0
\end{equation}
in view of $\rho_{m_k}\in L^\infty(0,T; {\cal
M}(b^{-s}_{q,\infty}(\Om)))$ and  $u_{m_kt},
(u_{m_k}\cdot\na)u_{m_k}, \na p_{m_k}\in L^\infty_{s/2}(0,T;
b^{-s}_{q,\infty}(\Om))$ by Lemma \ref{L4.3}. In view of the fact
that each term of \eq{4.14} belongs to $L^1(0,T;
b^{-s}_{q,\infty}(\Om))$ by \eq{4.25} and
$(b^{-s}_{q,\infty}(\Om))'=B^s_{q',1}(\Om)$, we get by testing
\eq{4.14} with arbitrary $\vp\in C_0^\infty([0,T)\ti\Om)^n$, $\div
\vp=0$, that
$$
\begin{array}{l}
\int_0^T\lan \rho_{m_k} u_{m_kt}-\mu \Da u_{m_k}+\rho_{m_k}
(u_{m_k}\cdot\na) u_{m_k},\vp\ran_{b^{-s}_{q,\infty}(\Om),
B^s_{q',1}(\Om)}\,dt=0.
\end{array}
$$
By Lemma \ref{L4.3} we have
$$u_m\in L^\infty_{1-\tau+s/2}(0,T; B^{2\tau-s}_{q,1}(\Om)),\;\forall m\in\N,$$
with $\tau=s-\theta+1/2$ provided $s-1/2<\theta<(s+1)/2$. Hence, if
$(1+s)/2-\theta)p'<1$, i.e., if $p>2/(2\theta-s+2)$, then
 \begin{equation}
 \label{E4.23}
u_{m_k}\in L^{p'}(0,T;H^{1-2\theta+s}_{q,0}(\Om))\subset
L^{p'}(0,T;H^{1-2\theta+s}_{q',0}(\Om)).
 \end{equation}
  On the other hand, by
\eq{4.13}, we have
$$\rho_{m_kt}\vp\in L^p(0,T;H^{-1+2\theta-s}_{q}(\Om))$$
if $p\in (1, 2/(2\theta-s))$, $\theta\in (s/2, s/2+1/(2q))$.
Therefore, if
 \begin{equation}
 \label{E4.22}
p\in \big(\frac{2}{2\theta-s+2}, \frac{2}{2\theta-s}\big),\quad
 \theta\in \big(\max\{\frac{s}{2}, s-\frac{1}{2}\}, \frac{s}{2}+\frac{1}{2q}\big),
 \end{equation}
then we have
$$\begin{array}{rcl}
\int_0^T\lan\rho_{m_k} u_{m_kt}, \vp\ran_{b^{-s}_{q,\infty}(\Om),
B^s_{q',1}(\Om)}\,dt &=& -\lan u_{m_k},\rho_{m_kt}\vp
\ran_{L^{p'}(0,T;H^{1-2\theta+s}_{q',0}(\Om)),L^p(0,T;H^{-1+2\theta-s}_{q}(\Om))}\ek
&&-\int_0^T\int_\Om\rho_{m_k}u_{m_k}\cdot\vp_t\,dxdt -\int_\Om
\rho_0u_0\cdot\vp(0,\cdot)\,dx.
\end{array}$$
 Moreover, for $p$ and $\theta$ satisfying \eq{4.22} we have
$$\begin{array}{l}
\int_0^T \lan \rho_{m_k} (u_{m_k}\cdot\na)
u_{m_k},\vp\ran_{b^{-s}_{q,\infty}(\Om), B^s_{q',1}(\Om)}\,dt
=-\int_0^T\int_\Om \rho_{m_k} u_{m_k}\otimes
u_{m_k}\cdot\na\vp\,dxdt\ek\hspace{1cm} -\lan u_{m_k},
(u_{m_k}\cdot\na)
\rho_{m_k}\vp\ran_{L^{p'}(0,T;H^{1-2\theta+s}_{q',0}(\Om)),L^p(0,T;H^{-1+2\theta-s}_{q}(\Om))},
\end{array}$$
in view of \eq{4.23} and the fact that $(u_{m}\cdot\na) \rho_{m}\in
L^p(0,T;H^{-1+2\theta-s}_{q}(\Om))$ holds true since
$(u_{m}\cdot\na) \rho_{m}=\div(\rho_mu_m)$ and $\rho_mu_m\in
L^p(0,T;H^{2\theta-s}_{q}(\Om))$ due to $\rho_m\in L^\infty(0,T;
{\cal M}(B^{2\theta-s}_{q,1}(\Om))$, $u_m\in
L^p(0,T;B^{2\theta-s}_{q,1}(\Om))$ by Lemma \ref{L4.3}.

 Therefore, we have
\begin{equation}
\label{E4.17}
\begin{array}{l}
0=\int_0^T\int_{\Om}[\rho_{m_k} u_{m_k}\cdot\vp_t+\mu u_{m_k}\cdot
\Da\vp +\rho_{m_k} u_{m_k}\otimes u_{m_k}\cdot\na\vp]\,dxdt
+\int_\Om \rho_0u_0\cdot\vp(0,\cdot)\,dx\ek\hspace{1cm} +\lan
u_{m_k},(\rho_{m_kt}+(u_{m_k}\cdot\na)\rho_{m_k})\vp
\ran_{L^{p'}(0,T;H^{1-2\theta+s}_{q',0}(\Om)),L^p(0,T;H^{-1+2\theta-s}_{q}(\Om))}\ek
=\int_0^T\int_{\Om}[\rho_{m_k} u_{m_k}\cdot\vp_t+\mu u_{m_k}\cdot
\Da\vp +\rho_{m_k} u_{m_k}\otimes u_{m_k}\cdot\na\vp]\,dxdt
+\int_\Om \rho_0u_0\cdot\vp(0,\cdot)\,dx\ek\hspace{1cm} +\lan
u_{m_k},((u_{m_k}-u^{(m_k-1)})\cdot\na)\rho_{m_k}\vp
\ran_{L^{p'}(0,T;H^{1-2\theta+s}_{q',0}(\Om)),L^p(0,T;H^{-1+2\theta-s}_{q}(\Om))},\ek
\hfill\forall \vp\in C_0^\infty([0,T)\ti\Om)^n\,\,(\div \vp=0).
\end{array}
\end{equation}
\par\bigskip
For the estimate of the last term on the right-hand side of
\eq{4.17} we need the following lemma.
\begin{lemma}
\label{L4.6} {\rm Let $p$ and $\theta$ satisfy \eq{4.22} and put
$$R_k:=\lan u_{m_k},((u_{m_k}-u^{(m_k-1)})\cdot\na)\rho_{m_k}\vp
\ran_{L^{p'}(0,T;H^{1-2\theta+s}_{q',0}(\Om)),L^p(0,T;H^{-1+2\theta-s}_{q}(\Om))},$$
where $\vp\in C_0^\infty([0,T)\ti\Om)^n$ with $\div \vp=0$. Then
$R_k$ tends to $0$ as $k\ra\infty$.
 }
\end{lemma}
{\bf Proof:} Let $Q_T=(0,T)\ti\Om$.
In view of Remark \ref{R4.4} and Lemma \ref{L4.3} we get that for any sufficiently
small $\da>0$
$$\begin{array}{rl}
|R_k|&=|\lan \rho_{m_k}(u_{m_k}-u^{(m_k-1)}), \na(u_{m_k}\cdot\vp)\ran_{Q_T}|\ek
&\leq \|\rho_{m_k}(u_{m_k}-u^{(m_k-1)})\|_{L^{q(n+2)/n-\da}(Q_T)}
\|\na(u_{m_k}\cdot\vp)\|_{(L^{q(n+2)/n-\da}(Q_T))'}\ek
&\leq \|\rho_0\|_{L^\infty(\Om)}\|u_{m_k}-u^{(m_k-1)}\|_{L^{q(n+2)/n-\da}(Q_T)}
\|\na (u_{m_k}\cdot\vp)\|_{(L^{q(n+2)/n-\da}(Q_T))'}.
\end{array}$$
Then, since $\|\na (u_{m_k}\cdot\vp)\|_{(L^{q(n+2)/n-\da}(Q_T))'}$
is bounded with respect to $k\in\N$, see \eq{4.11}, the proof of the
lemma is complete if we show that
\begin{equation}
\label{E4.18} \|u_{m_k}-u^{(m_k-1)}\|_{L^{q(n+2)/n-\da}(Q_T)}\ra
0\quad (k\ra \infty).
\end{equation}
Note that
$$\|u_{m_k}-u^{(m_k-1)}\|_{L^{q(n+2)/n-\da}(Q_T)}\leq \|u_{m_k}-u\|_{L^{q(n+2)/n-\da}(Q_T)}
+\|u^{(m_k-1)}-u\|_{L^{q(n+2)/n-\da}(Q_T)},$$ where the first term
on the right-hand side tends to $0$ as $k\ra\infty$ due to
\eq{4.20}.

In order to show
\begin{equation}
\label{E4.27} \|u^{(m_k-1)}-u\|_{L^{q(n+2)/n-\da}(Q_T)}\ra 0\quad
(k\ra\infty),
\end{equation}
 we write
$$u^{(m_k-1)}-u=\eta_{m_k-1}\star (\bar{u}_{m_k-1}-\tilde{u})|_\Om
+(\eta_{m_k-1}\star\tilde{u}-\tilde{u})|_\Om.$$
Here, $L^{q(n+2)/n-\da}(Q_T)$-norm of $(\eta_{m_k-1}\star\tilde{u}-\tilde{u})|_\Om$
 obviously tends to $0$ as $k\ra\infty$.
 The $L^{q(n+2)/n-\da}(Q_T)$-norm of $\eta_{m_k-1}\star (\bar{u}_{m_k-1}-\tilde{u})|_\Om$
 goes to zero as $k\ra\infty$ since
$\|\eta_{m_k}\star\cdot\|_{{\cal L}(L^{q(n+2)/n-\da}(Q_T))}$ is uniformly bounded with
respect to $k\in\N$ and
$$\begin{array}{l}
\|\bar{u}_{m_k}-\tilde{u}\|_{L^{q(n+2)/n-\da}(Q_T)}\ek\hspace{1cm}
\leq \|u_{m_k}-u\|_{L^{q(n+2)/n-\da}((0,T)\ti \Om'_{m_k})}
+\|u\|_{L^{q(n+2)/n-\da}((0,T)\ti (\Om\setminus\Om'_{m_k}))}\ra 0
\end{array}$$
as $k\ra\infty$, where $\Om'_{m_k}=\{x\in
\Om;\text{dist}(x,\pa\Om)\geq\frac{1}{m_k}\}$. Thus, \eq{4.27} and,
consequently, \eq{4.18} are proved. The proof of the lemma comes to
end.\qed
\par\bigskip

Let us continue the proof of existence part of Theorem \ref{T1.2}.
In \eq{4.17}, we get easily that for any sufficiently small $\da>0$
$$\int_0^T\int_{\Om}(\rho_{m_k} u_{m_k}\cdot\vp_t+\mu u_{m_k}\cdot \Da\vp)\,dxdt
\ra \int_0^T\int_{\Om}(\rho u\cdot\vp_t+\mu u\cdot \Da\vp)\,dxdt$$
as $k\ra\infty$ due to  $*$-weak convergence $\rho_{m_k}
\rightharpoonup \rho$  in
           $L^\infty(Q_T)$
and strong convergence $u_{m_k} \ra u$ in $L^{q(n+2)/n-\da}(Q_T)$,
see \eq{4.8} and \eq{4.20}. Note that $q(n+2)/n-\da\geq 2$ and hence
$u_{m_k}\otimes u_{m_k}\ra u\otimes u$ in $L^1(\Om)$ as
$k\ra\infty$.
 Therefore,
$$\int_0^T\int_{\Om}
\rho_{m_k} (u_{m_k}\otimes u_{m_k})\cdot\na\vp\,dxdt \ra
\int_0^T\int_{\Om}
\rho (u\otimes u)\cdot\na\vp\,dxdt\quad\text{as }k\ra\infty.$$

Thus, we get \eq{1.3}  letting $k\ra\infty$ in \eq{4.17}.
\par\bigskip
 Next, in order to show \eq{1.2},
 test the first equation of \eq{4.4}
 with $\psi\in C_0^1([0,T)\ti\Om)$ to get
\begin{equation}
\label{E4.16}
\int_0^T\int_{\Om}(\rho_{m_k}\psi_t+\rho_{m_k} u^{(m_k-1)}\cdot \na\psi)\,dxdt
+\int_\Om \rho_0\psi(0,\cdot)\,dx=0,\quad \forall\psi\in C_0^1([0,T)\ti\Om).
\end{equation}
Obviously,
$$\int_{Q_T} \rho_{m_k}\psi_t\,dxdt\ra \int_{Q_T} \rho\psi_t\,dxdt\quad (k\ra\infty).$$
Moreover, for all $\psi\in C_0^1([0,T)\ti\Om)$ we get in view of \eq{4.27} and \eq{4.8}
that
$$\begin{array}{l}
\big|\int_{Q_T}(\rho_{m_k} u^{(m_k-1)}-\rho u)\cdot\na\psi\,dxdt\big|\ek
 \leq \big|\int_{Q_T}(u^{(m_k-1)}- u)\cdot(\rho_{m_k} \na\psi)\,dxdt\big|
 + \big|\int_{Q_T}(\rho_{m_k}-\rho)u\cdot\na\psi\,dxdt\big|\ra 0
\end{array}$$
as $k\ra\infty$.
Thus, $(\rho,u)$ satisfies \eq{1.2} in the limiting case $k\ra\infty$ in \eq{4.16}.

Finally, the proof of existence part for Theorem \ref{T1.2} is
completed. \qed
\par\bigskip
\subsection{Proof of uniqueness}
For the proof of uniqueness part for Theorem \ref{T1.2} we need the
following statement.
\begin{lemma}
\label{L4.7}
 {\rm
 Let $\Om$ be a Lipschitz domain of $\R^n$, $n\in\N$, and
let $q>n$ and $s\in (0, \frac{q-n}{2q-n})$. Then there hold the
followings:
\begin{itemize}
\item[(i)] It holds
$$ B^{1-s}_{q,\infty}(\Om)\cdot B^s_{q',1}(\Om)\hookrightarrow B^s_{q',1}(\Om), q'=q/(q-1),$$
and
$$\|f\vp \|_{B^s_{q',1}(\Om)}\leq c(q,s,\Om)\|f\|_{B^{1-s}_{q,\infty}(\Om)}
\|\vp\|_{B^s_{q',1}(\Om)},\quad \forall f\in
B^{1-s}_{q,\infty}(\Om), \vp\in B^s_{q',1}(\Om).$$

\item[(ii)] If $s\in (0, \frac{1}{q'})$, then
$$ b^{1-s}_{q,\infty}(\Om)\cdot b^{-s}_{q,\infty}(\Om)\hookrightarrow b^{-s}_{q,\infty}(\Om)$$
and
$$\|f g\|_{b^{-s}_{q,\infty}(\Om)}\leq
k_0(q,s,\Om)\|f\|_{b^{1-s}_{q,\infty}(\Om)}\|g\|_{b^{-s}_{q,\infty}(\Om)},
\quad \forall f\in b^{1-s}_{q,\infty}(\Om), g\in
b^{-s}_{q,\infty}(\Om).$$

\item[(iii)]  If $s\in (0, \frac{1}{q'})$, then
$$ b^{1-s}_{q,\infty}(\Om)\cdot b^{1-s}_{q,\infty}(\Om)\hookrightarrow b^{1-s}_{q,\infty}(\Om)$$
and
$$\|f g\|_{b^{1-s}_{q,\infty}(\Om)}\leq
k_1(q,s,\Om)\|f\|_{b^{1-s}_{q,\infty}(\Om)}\|g\|_{b^{1-s}_{q,\infty}(\Om)},
\quad \forall f,g\in b^{1-s}_{q,\infty}(\Om).$$
\end{itemize}
 }
\end{lemma}
{\bf Proof:} -- {\it Proof of (i)}: Thanks to $q>n$, it is clear
that
\begin{equation}
\label{E4.39}
H^1_q(\Om)\cdot H^1_{q'}(\Om)\hookrightarrow H^1_{q'}(\Om).
\end{equation}
Note that $\frac{n}{q}<\frac{1-2s}{1-s}$ due to the assumption on $s$.
Let us choose $\al\in (\frac{n}{q}, \frac{1-2s}{1-s})$.
Then, $H^\al_q(\Om)\hookrightarrow L^\infty(\Om)$ and
\begin{equation}
\label{E4.40}
H^\al_q(\Om)\cdot L^{q'}(\Om)\hookrightarrow L^{q'}(\Om).
\end{equation}
From \eq{4.39} and \eq{4.40} we get by bilinear interpolation that
$$(H^\al_q(\Om), H^1_q(\Om))_{s,1}\cdot (L^{q'}(\Om), H^1_{q'}(\Om))_{s,1}
\hookrightarrow (L^{q'}(\Om), H^1_{q'}(\Om))_{s,1},
$$
i.e.,
$$
B^{(1-s)\al+s}_{q,1}(\Om)\cdot B^s_{q',1}(\Om)\hookrightarrow B^s_{q',1}(\Om).
$$
By the way, we have $B^{1-s}_{q,\infty}(\Om) \hookrightarrow
B^{(1-s)\al+s}_{q,1}(\Om)$ since $\Om$ is bounded and
$1-s>(1-s)\al+s$ due to $\al<\frac{1-2s}{1-s}$. Hence we get the
conclusion.

-{\it Proof of (ii)}: Let $f\in B^{1-s}_{q,\infty}(\Om)$, $g\in
b^{-s}_{q,\infty}(\Om)$. Recall that
$(b^{-s}_{q,\infty}(\Om))'=B^s_{q',1}(\Om)$ and
$(B^s_{q',1}(\Om))'=B^{-s}_{q,\infty}(\Om)$ for $0<s<1/q'$. Then, by
the first assertion of the lemma we have
$$\begin{array}{rcl}
|\lan fg, \vp\ran_{\Om}|
&=&|\lan g, f \vp\ran_{b^{-s}_{q,\infty}(\Om), (b^{-s}_{q,\infty}(\Om))'}|\ek
&\leq& \|g\|_{b^{-s}_{q,\infty}(\Om)}\|f\vp\|_{B^{s}_{q',1}(\Om)}\ek
&\leq &c(q,s,\Om)\|g\|_{b^{-s}_{q,\infty}(\Om)}\|f\|_{B^{1-s}_{q,\infty}(\Om)}
\|\vp\|_{B^{s}_{q',1}(\Om)}
\end{array}$$
for any $\vp\in B^{s}_{q',1}(\Om)$. Hence we have $fg\in
B^{-s}_{q,\infty}(\Om)$ and
\begin{equation}
\label{E4.41} \|f g\|_{B^{-s}_{q,\infty}(\Om)} \leq
k_0\|f\|_{B^{1-s}_{q,\infty}(\Om)}\|g\|_{b^{-s}_{q,\infty}(\Om)}
\end{equation}
with some $k_0=k_0(q,s,\Om)>0$.

Thus it remained to prove $fg\in b^{-s}_{q,\infty}(\Om)$ provided
$f\in b^{1-s}_{q,\infty}(\Om)$. Let $\{f_m\}, \{g_m\}\subset
C^\infty(\bar\Om)$ be sequences converging to $f$ and $g$ in
$b^{1-s}_{q,\infty}(\Om)$ and $b^{-s}_{q,\infty}(\Om)$,
respectively. Then, by \eq{4.41} we have
$$
\|f_m g_m-f g\|_{B^{-s}_{q,\infty}(\Om)}
\leq c(\|f_m-f\|_{B^{1-s}_{q,\infty}(\Om)}\|g_m\|_{b^{-s}_{q,\infty}(\Om)}+
\|f\|_{B^{1-s}_{q,\infty}(\Om)}\|g_m-g\|_{b^{-s}_{q,\infty}(\Om)})\ra 0
$$
as $m\ra \infty$. Hence, $fg\in b^{-s}_{q,\infty}(\Om)$ and the
assertion follows.

--{\it Proof of (iii)}: By a similar argument to prove (ii) one can
show that, if $s\in (0, \frac{1}{q'})$,
$$ H^{1-s}_{q}(\R^n)\cdot H^{-s}_{q}(\R^n)\hookrightarrow H^{-s}_{q}(\R^n)$$
and
$$\|f g\|_{H^{-s}_{q}(\R^n)}\leq c\|f\|_{H^{1-s}_{q}(\R^n)}\|g\|_{H^{-s}_{q}(\R^n)},\,\,
\forall f\in H^{1-s}_{q}(\R^n), g\in f\in H^{-s}_{q}(\R^n).$$
Therefore, in view of the fact that
$$\begin{array}{l}
H^{1-s}_{q}(\R^n)=\{u\in H^{-s}_{q}(\R^n): \na u\in
H^{-s}_{q}(\R^n)\},\ek
 \|u\|_{H^{1-s}_{q}(\R^n)}\sim \|u\|_{H^{-s}_{q}(\R^n)}+\|\na
 u\|_{H^{-s}_{q}(\R^n)},
 \end{array}
$$
we have
\begin{equation}
\label{E5.27}
\begin{array}{l}
H^{1-s}_{q}(\R^n)\cdot H^{1-s}_{q}(\R^n)\hookrightarrow
H^{1-s}_{q}(\R^n),\ek
 \|f g\|_{H^{1-s}_{q}(\R^n)}\leq
c\|f\|_{H^{1-s}_{q}(\R^n)}\|g\|_{H^{1-s}_{q}(\R^n)},\,\, \forall
f,g\in H^{1-s}_{q}(\R^n).
\end{array}
\end{equation}
Then, \eq{5.27} for general Lipschitz domains can be proved easily
by the Sobolev extension theorem. Finally, the assertion follows by
bilinear interpolation.
 \qed
\par\bigskip
{\bf Proof of Theorem \ref{T1.2}: uniqueness part}
\par\bigskip
The uniqueness proof relies on a Lagrangian coordinates approach,
 using, in principle, the same idea as that on pages 29-31 of \cite{DaZh13}
 but being based on pointwise multiplier in little Nicoskii spaces.

First let us recall some facts concerning Lagrangian coordinates.
Let $u$ be a velocity vector field such that
\begin{equation}
\label{E4.32}
u\in L^1_{loc}(0,T; C_{Lip}(\R^n))
\end{equation}
and let $X(t,y)$ be the solution to the
ordinary differential system on $(0,T)$:
\begin{equation}
\label{E4.28}
\frac{dX}{dt}=u(t,X)\quad t\in (0,T), \quad X(0)=y\in \R^n.
\end{equation}
Then, Eulerian coordinates $x$ and Lagrangian coordinates $y$
are related by
\begin{equation}
\label{E4.29} x=X(t,y)=y+\int_0^t u(\tau,X(\tau,y))\,d\tau.
\end{equation}
Thus, given a vector field $u$ satisfying \eq{4.32}, a unique  $C^1$-semiflow $X$
is determined. In particular, if $u|_{\pa\Om}=0$
for $\Om\subset\R^n$ with $\pa\Om\in C^1$ and $y\in\Om$, then
$X(t,y)\in\Om$ for all $t>0$. Therefore, if
\begin{equation}
\label{E4.33} u\in L^1_{loc}(0,T; C_{Lip}(\Om)),\quad u|_{\pa\Om}=0
\end{equation}
instead of \eq{4.32}, then a unique $C^1$-semiflow $X$ in $\Om$
given by \eq{4.29} is generated. Note that $W^{1,\infty}(\Om)\subset
C_{Lip}(\Om)$. Let $Y(t,\cdot)$ be the inverse mapping of
$X(t,\cdot)$, then
$$D_xY(t,x)=(D_yX(t,y))^{-1}=:A(t,y),$$
where and in what follows we use the notation
$(\na u)_{i,j}=(\pa_iu^j)_{1\leq i,j\leq n}, Du=(\na u)^T$.
Let
$v(t,y):=u(t, X(t,y))$.
Then,
\begin{equation}
\label{E4.30}
\na_x u(t,x)=A^T \na_y v(t,y),\quad \text{div}_x u(x,t)=\text{div}_y(A v(t,y)),
\end{equation}
see \cite{DaMu12}, Appendix.
In view of \eq{4.30}, we use the notation
$$\na_u w:= A^T \na_y w,\quad \text{div}_u w:= \text{div}_y(Aw ),
\quad\Da_uw:=\text{div}_u(\na_u w).$$

\par\bigskip
Now, let $\{\rho, u\}$ be a solution to \eq{1.1} in the sense of
Definition \ref{D1.1} with $q>n$ and $s\in (0,
\min\{\frac{n}{q},1-\frac{1}{q},\frac{q-n}{2q-n}\})$, whose
existence is guaranteed by Theorem \ref{T1.2}. Then, $u\in L^1(0,T;
W^{1,\infty}(\Om))$
 due to $u\in L^\infty_{s/2}(0,T; b^{2-s}_{q,\infty}(\Om))$
and $1-s-n/q>0$.

Let $a(t,x):=\frac{\bar\rho}{\rho(t,x)}-1$ and
$$b(t,y):=a(t,X(t,y)),\quad v(t,y):=u(t,X(t,y))\quad\text{and}\quad Q(t,y):=P(t,X(t,y)).$$
Then
we get from \eq{1.2} that
$b_t=0$ in the sense of distribution.
In fact, given any $\tilde\psi\in C^1_0(\Om\ti (0,T))$, for
$\psi(t,x):=\tilde\psi(t,Y(t,x))$ we have $\tilde\psi(t,y)=\psi(t, X(t,y))$ and $\psi\in
C^1_0(\Om\ti(0,T))$. Hence, we have
$$\begin{array}{rcl}
\int_{\Om\ti(0,T)}b\tilde\psi_t\,dydt&=&
\int_{\Om\ti(0,T)}a(t, X(t,y)) \frac{\pa\psi}{\pa t}(t,X(t,y))\,dydt\ek
&=&\int_{\Om\ti(0,T)}a(t, X(t,y)) (\psi_t+u\cdot\na\psi)(t,X(t,y))\,dydt\ek
&=&\int_{\Om\ti(0,T)}a(t,x) (\psi_t+u\cdot\na\psi)(t,x)|D_y X(t,y)|\,dxdt\ek
&=&\int_{\Om\ti(0,T)}a(t,x) (\psi_t+u\cdot\na\psi)(t,x)\,dxdt=0,
\end{array}$$
where we used that the Jacobian $|D_y X(t,y)|=1$ due to the
solenoidal condition for the vector filed $u$, i.e., $\div u=0$.
Therefore, $b(t,y)\equiv b(0,y)\equiv a_0(y)$ for each $y\in \Om$.

Assuming $\mu=1$ without loss of generality, $\{v,Q\}$ solves the
system:
\begin{equation}
\label{E4.34}
\begin{array}{rcl}
     v_t-(1+b)(\Da_u v -\na_u Q)  = 0 \,\, &\text{in}&(0,T)\ti\Om,\ek
      \text{div}_u v = 0 \,\, &\text{in} &(0,T)\ti\Om,\ek
     v= 0\, &\text{on}&(0,T)\ti\pa\Om,\ek
     v(0,y)=u_0&\text{in}&\Om.
     \end{array}
\end{equation}

Now, let $\{\rho_i, u_i, P_i\}$, $i=1,2$,  be two solutions to
\eq{1.1}, and let $X_i$ be the semiflow of $u_i$ and let $(b_i, v_i,
Q_i)$ be the corresponding density perturbation, velocity  and
pressure in the Lagrangian coordinates. Then, for $\da v=v_1-v_2,
\da Q=Q_1-Q_2$ one has
\begin{equation}
\label{E4.35}
\begin{array}{rclll}
     (\da v)_t-\Da \da v +\na \da Q  &=& a_0(\Da\da v-\na\da Q)+\da F \,\,
      &\text{in}&(0,T)\ti\Om,\ek
      \div \da v &=& \da g=\div \da R \,\, &\text{in} &(0,T)\ti\Om,\ek
     \da v&=& 0\, &\text{on}&(0,T)\ti\pa\Om,\ek
     \da v(0,y)&=& 0& \text{in}&\Om,
\end{array}
\end{equation}
where
\begin{equation}
\label{E4.46}
\begin{array}{l}
\da F=\da f_1+ \da f_2,\ek
 \da f_1:= (1+a_0)[(\text{Id}-A_2^T)\na\da
Q-\da A\na Q_1]\quad\text{with}\quad \da A:= A_2-A_1,
 \ek
  \da f_2:=(1+a_0)\div [(A_2^T A_2-\Id)\na\da v-(A_2^T A_2-A_1^T A_1)\na
  v_1],
\ek \da g:= (\Id-A_2^T): \na \da v-(\da A)^T:\na v_1,
 \ek \da R:=(\Id-A_2)\da v-\da A v_1.
\end{array}
\end{equation}
 Note that $\da R(0)=0$. Therefore, by
Theorem \ref{T2.7} we have
\begin{equation}
\label{E4.36}
\begin{array}{l}
\|\da v_t,\na^2 \da v, \na \da
Q\|_{L^\infty_{s/2}(0,t;b^{-s}_{q,\infty}(\Om))} +\|\da
v\|_{L^{\infty}(0,t; b^{0}_{q,\infty}(\Om))}\ek\hspace{1cm} \leq K
\|a_0(\Da\da v-\na\da Q),\da F, \na \da g, (\da
R)_t\|_{L^\infty_{s/2}(0,t;b^{-s}_{q,\infty}(\Om))}, \forall t>0,
\end{array}
\end{equation}
with constant $K>0$ independent of $t$.

From now on, let us get estimate of the right-hand side of
\eq{4.36}.

Thanks to $\na v_i\in L^\infty_{s/2}(0,T; b^{1-s}_{q,\infty}(\Om)),
i=1,2$, and  $b^{1-s}_{q,\infty}(\Om)\hookrightarrow L^\infty(\Om)$
there is some $m_0>0$ such that, if $|t_2-t_1|<m_0$,
\begin{equation}
\label{E4.44}
 \big\|\int_{t_1}^{t_2} \na v_i(t,y)\,dt\big\|_{L^\infty(\Om)}<1,\,\,
   \int_{t_1}^{t_2} \|\na
v_i(t,y)\|_{b^{1-s}_{q,\infty}(\Om)}\,dt<\frac{1}{2k_2},i=1,2,
\end{equation}
where $k_2:=\max\{k_0(q,s,\Om), k_1(q,s,\Om)\}$ with constants
$k_0$, $k_1$ appearing in Lemma \ref{L4.7}. Throughout the proof, we
assume that $0<t<m_0$. Then we have
\begin{equation}
\label{E4.21}
 A_i(t,y)=(\Id+C_i(t))^{-1}=\sum_{k\geq 0}(-1)^k
C_i(t)^k, i=1,2,
\end{equation}
where $C_i(t):=\int_0^t D v_i(\tau)\,d\tau, i=1,2$. Hence one has
\begin{equation}
\label{E4.45} \da A(t)=h_1(t)\int_0^t D\da v (\tau)\,d\tau,t\in
(0,m_0),
\end{equation}
with
\begin{equation}
\label{E4.62} h_1(t):=\sum_{k\geq
1}(-1)^k\sum_{j=0}^{k-1}C_1(t)^jC_2(t)^{k-1-j},\,\,t \in (0,m_0).
\end{equation}
Since
\begin{equation}
\label{E4.59}
\bar{C}_i(t):=\|C_i\|_{L^\infty(0,t;b^{1-s}_{q,\infty}(\Om))}<\frac{1}{2k_2(q,s,\Om)},
\end{equation}
we get by Lemma \ref{L4.7} (ii) that
\begin{equation}
\label{E4.63n} \|h_1\|_{L^\infty(0,t;{\cal
M}(b^{-s}_{q,\infty}(\Om)))}\leq \sum_{k\geq 1}
k_0(q,s,\Om)^{k-1}\sum_{j=0}^{k-1}\bar{C}_1(t)^j\bar{C}_2(t)^{k-1-j}<C
\end{equation}
with $C=C(q,s,\Om)$ and, in particular,
\begin{equation}
  \label{E4.60}
 \|h_1\|_{L^\infty(0,t;{\cal M}(b^{-s}_{q,\infty}(\Om)))}\leq
 c(q,s,\Om)t^{s/2}\|v_1,
 v_2\|_{L^\infty_{s/2}(0,t;b^{2-s}_{q,\infty}(\Om))}.
 \end{equation}
Hence,
\begin{equation}
\label{E4.58} \|\da A\|_{L^\infty(0,t;b^{-s}_{q,\infty}(\Om))}\leq
 ct^{s/2}\|\da v\|_{L^\infty_{s/2}(0,t; b^{1-s}_{q,\infty}(\Om))}
\end{equation}
with $c=c(q,s,\Om)>0$. In the same way, using Lemma \ref{L4.7} (iii)
we can obtain
\begin{equation}
\label{E4.61}
 \|\da A\|_{L^\infty(0,t;b^{1-s}_{q,\infty}(\Om))}\leq
 ct^{s/2}\|\da v\|_{L^\infty_{s/2}(0,t; b^{2-s}_{q,\infty}(\Om))}
\end{equation}
with $c=c(q,s,\Om)>0$.

Since $a_0(y)=\frac{\bar\rho-\rho_0(y)}{\rho_0(y)}$, there is some
$\da_1=\da_1(\Om)>0$ such that if
\begin{equation}
\label{E4.51}
\frac{\rho_{02}-\rho_{01}}{\rho_{01}}<\da_1,
\end{equation}
 then
\begin{equation}
\label{E4.37} \|a_0(\Da \da v-\na\da
Q)\|_{L^\infty_{s/2}(0,t;b^{-s}_{q,\infty}(\Om))} \leq
\frac{1}{2K}\|\na^2 \da v, \na \da
Q\|_{L^\infty_{s/2}(0,t;b^{-s}_{q,\infty}(\Om))}.
\end{equation}

Next, let us get estimate of $\|\na\da
g\|_{L^\infty_{s/2}(0,t;b^{-s}_{q,\infty}(\Om))}$. From \eq{4.46}
one has
\begin{equation}
\label{E4.38}
\begin{array}{l}
\|\na\da g\|_{L^\infty_{s/2}(0,t;b^{-s}_{q,\infty}(\Om))} \ek \leq
\|\na A^T_2\otimes \na \da v, (\Id-A^T_2) \otimes \na^2 \da v,
\na\da A^T\otimes \na v_1, \da A^T\otimes\na^2
v_1\|_{L^\infty_{s/2}(0,t;b^{-s}_{q,\infty}(\Om))},
\end{array}
\end{equation}
where the right-hand side can be estimated as below.

 Since $DA_i=\sum_{k\geq 1}(-1)^kkC_i^{k-1}DC_i$ due to \eq{4.21},
 we get by Lemma \ref{L4.7} (ii) and \eq{4.59} that
\begin{equation}
\label{E4.53} \begin{array}{rcl}
\|\na A^T_2\|_{L^\infty(0,t;
b^{-s}_{q,\infty}(\Om))}
&\leq & c\|\na^2 v_2\|_{L^1(0,t; b^{-s}_{q,\infty}(\Om))} \ek
 &\leq & ct^{s/2}\|\na^2 v_2\|_{L^\infty_{s/2}(0,t;
b^{-s}_{q,\infty}(\Om))}
\end{array}
\end{equation}
and
\begin{equation}
\label{E4.42}
\begin{array}{rl}
\|\na A^T_2\otimes \na \da
v\|_{L^\infty_{s/2}(0,t;b^{-s}_{q,\infty}(\Om))} &\leq c\|\na
A^T_2\|_{L^\infty(0,t;b^{-s}_{q,\infty}(\Om))}
 \|\na \da v\|_{L^\infty_{s/2}(0,t;b^{1-s}_{q,\infty}(\Om))}
\ek
  &\leq
ct^{s/2}\|\na^2 v_2\|_{L^\infty_{s/2}(0,t; b^{-s}_{q,\infty}(\Om))}
  \|\da v\|_{L^\infty_{s/2}(0,t;b^{2-s}_{q,\infty}(\Om))}
\end{array}
\end{equation}
with $c=c(q,s,\Om)$. Furthermore, by Lemma \ref{L4.7} (ii) and
\eq{4.59} we have
\begin{equation}
\label{E4.54}
\begin{array}{rcl}
 \|\Id-A_i\|_{L^\infty(0,t;b^{1-s}_{q,\infty}(\Om))}
 &= &\|\sum_{k\geq
1}(-1)^kC_i(t)^k\|_{L^\infty(0,t;b^{1-s}_{q,\infty}(\Om))}\ek
 &\leq &
c\|\na v_i\|_{L^1(0,t;b^{1-s}_{q,\infty}(\Om))}, i=1,2,
 \end{array}
\end{equation}
and, consequently,
\begin{equation}
\label{E4.43}
\begin{array}{l}
\|(\Id-A^T_2) \otimes \na^2 \da
v\|_{L^\infty_{s/2}(0,t;b^{-s}_{q,\infty}(\Om))}\ek\hspace{1cm} \leq
c\|\Id-A^T_2\|_{L^\infty(0,t;b^{1-s}_{q,\infty}(\Om))} \|\na^2 \da
v\|_{L^\infty_{s/2}(0,t;b^{-s}_{q,\infty}(\Om))}\ek\hspace{1cm} \leq
c\|\na v_2\|_{L^1(0,t;b^{1-s}_{q,\infty}(\Om))} \|\na^2 \da
v\|_{L^\infty_{s/2}(0,t;b^{-s}_{q,\infty}(\Om))}\ek\hspace{1cm} \leq
ct^{s/2}\|v_2\|_{L^\infty_{s/2}(0,t; b^{2-s}_{q,\infty}(\Om))}
  \|\da v\|_{L^\infty_{s/2}(0,t;b^{2-s}_{q,\infty}(\Om))}
\end{array}
\end{equation}
with $c=c(q,s,\Om)$.

Note that
$$\begin{array}{l} D\da A(t)= h_1(t) \int_0^t D^2\da v
(\tau)\,d\tau + h_2(t)\int_0^t D\da v (\tau)\,d\tau,
\end{array}
$$
where
$$\begin{array}{l}
h_2:=Dh_1=
  \sum_{k\geq 1}(-1)^k\sum_{j=0}^{k-1}\big(jC_1^{j-1}C_2^{k-1-j}DC_1+
  (k-1-j)C_1^{j}C_2^{k-2-j}DC_2\big).
\end{array}
$$
We get, in view of \eq{4.44}, that
$$\begin{array}{l}
\|h_2\|_{L^\infty(0,t;b^{-s}_{q,\infty}(\Om))}<c\|\na^2 v_1,\na^2
v_2\|_{L^1(0,t;b^{-s}_{q,\infty}(\Om))}<C_1
\end{array}
$$
 with $C_1=C_1(q,s,\Om)$. Hence,
$$\begin{array}{rl} \|D\da A\|_{L^\infty(0,t;
b^{-s}_{q,\infty}(\Om))}&\leq \|h_1\|_{L^\infty(0,t;{\cal
M}(b^{-s}_{q,\infty}(\Om)))}
 \|\int_0^t
D^2\da v (\tau)\,d\tau\|_{L^\infty(0,t; b^{-s}_{q,\infty}(\Om))}
 \ek
&\hspace{1cm}+ c\|h_2\|_{L^\infty(0,t;
b^{-s}_{q,\infty}(\Om))}\|\int_0^t D\da v
(\tau)\,d\tau\|_{L^\infty(0,t; b^{1-s}_{q,\infty}(\Om))}\ek
 &\leq ct^{s/2}\|\da v\|_{L^\infty_{s/2}(0,t; b^{2-s}_{q,\infty}(\Om))}
\end{array}$$
with $c=c(q,s,\Om)$. Therefore, by  Lemma \ref{L4.7} (ii) we have
\begin{equation}
\label{E4.47}
\begin{array}{rcl}
\|\na\da A^T\otimes \na
v_1\|_{L^\infty_{s/2}(0,t;b^{-s}_{q,\infty}(\Om))}&\leq & c\|\na\da
A\|_{L^\infty(0,t;b^{-s}_{q,\infty}(\Om))} \|\na
v_1\|_{L^\infty_{s/2}(0,t;b^{1-s}_{q,\infty}(\Om))}\ek\hspace{1cm}
&\leq& ct^{s/2}\| v_1\|_{L^\infty_{s/2}(0,t;
b^{2-s}_{q,\infty}(\Om))}
  \|\da v\|_{L^\infty_{s/2}(0,t;b^{2-s}_{q,\infty}(\Om))}.
\end{array}
\end{equation}

By \eq{4.58} we have
\begin{equation}
\label{E4.48}
\|\da A^T\otimes\na^2 v_1\|_{L^\infty_{s/2}(0,t;b^{-s}_{q,\infty}(\Om))}
\leq ct^{s/2}\| v_1\|_{L^\infty_{s/2}(0,t; b^{2-s}_{q,\infty}(\Om))}
  \|\da v\|_{L^\infty_{s/2}(0,t;b^{1-s}_{q,\infty}(\Om))}.
\end{equation}

Thus, from \eq{4.38}, \eq{4.42}, \eq{4.43}, \eq{4.47} and \eq{4.48}
we get
\begin{equation}
\label{E4.49}
\|\na\da g\|_{L^\infty_{s/2}(0,t;b^{-s}_{q,\infty}(\Om))}
\leq \eta_1(t)\|\da v\|_{L^\infty_{s/2}(0,t;b^{2-s}_{q,\infty}(\Om))}
\end{equation}
with $\eta_1(t):=ct^{s/2}\| v_1,v_2\|_{L^\infty_{s/2}(0,t; b^{2-s}_{q,\infty}(\Om))}\ra 0$
as $t\ra+0$.

Following the same procedure as the derivation of \eq{4.49}, we can
get estimate
\begin{equation}
\label{E4.50} \|\da
f_1\|_{L^\infty_{s/2}(0,t;b^{-s}_{q,\infty}(\Om))} \leq
ct^{s/2}\|\na^2 v_2, \na Q_1 \|_{L^\infty_{s/2}(0,t;
b^{-s}_{q,\infty}(\Om))}
  \|\na^2\da v, \na\da Q\|_{L^\infty_{s/2}(0,t;b^{-s}_{q,\infty}(\Om))}
\end{equation}
under the smallness assumption \eq{4.51}; we omit the details here.

On the other hand, we have
\begin{equation}
\label{E4.57} \|\da
f_2\|_{L^\infty_{s/2}(0,t;b^{-s}_{q,\infty}(\Om))} \leq \|\div
[(A_2^T A_2-\Id)\na\da v-(A_2^T A_2-A_1^T A_1)\na
v_1]\|_{L^\infty_{s/2}(0,t;b^{-s}_{q,\infty}(\Om))}.
\end{equation}
Here,
$$\begin{array}{l}
\|\div [(A_2^T A_2-\Id)\na\da v]\|_{L^\infty_{s/2}(0,t;b^{-s}_{q,\infty}(\Om))}
\ek\leq
(\|\na(A_2^T A_2)\|_{L^\infty(0,t; b^{-s}_{q,\infty}(\Om))}
            +\|A_2^T A_2-\Id\|_{L^\infty(0,t; b^{1-s}_{q,\infty}(\Om))})
        \|\da v\|_{L^\infty_{s/2}(0,t;b^{2-s}_{q,\infty}(\Om))}.
\end{array}$$
Note that, by Lemma \ref{L4.7} (ii), (iii), \eq{4.59} and \eq{4.53},
$$\begin{array}{rcl}
\|\na(A_i^T A_i)\|_{L^\infty(0,t; b^{-s}_{q,\infty}(\Om))}
&\leq& 2 \|\na A_i\|_{L^\infty(0,t; b^{-s}_{q,\infty}(\Om))}
 \|A_i\|_{L^\infty(0,t; b^{1-s}_{q,\infty}(\Om))}\ek
&\leq& c t^s \| v_i\|^2_{L^\infty_{s/2}(0,t; b^{2-s}_{q,\infty}(\Om))}, i=1,2,
\end{array}
$$
and, by Lemma \ref{L4.7} (iii),
$$\begin{array}{rcl}
\|A_2^T A_2-\Id\|_{L^\infty(0,t; b^{1-s}_{q,\infty}(\Om))} &\leq&
c\|v_2\|_{L^\infty(0,t;
b^{2-s}_{q,\infty}(\Om))}(1+\|v_2\|_{L^\infty(0,t;
b^{2-s}_{q,\infty}(\Om))})\ek
 &\leq& ct^{s/2}\|v_2\|_{L^\infty_{s/2}(0,t;
b^{2-s}_{q,\infty}(\Om))}(1+\|v_2\|_{L^\infty(0,t;
b^{2-s}_{q,\infty}(\Om))})
\end{array}
$$
in view of $$A_2^T(t) A_2(t)-\Id= \int_0^t (Dv_2+Dv_2^T)\,d\tau+
\int_0^tDv_2^T\,d\tau\cdot \int_0^tDv_2\,d\tau.$$ Therefore,
\begin{equation}
\label{E4.56}
\begin{array}{l}
\|\div [(A_2^T A_2-\Id)\na\da v]\|_{L^\infty_{s/2}(0,t;b^{-s}_{q,\infty}(\Om))}
\ek\leq
c t^{s/2} \| v_1,v_2\|_{L^\infty_{s/2}(0,t; b^{2-s}_{q,\infty}(\Om))}
(\|v_1,v_2\|_{L^\infty_{s/2}(0,t; b^{2-s}_{q,\infty}(\Om))}+1)
\|\da v\|_{L^\infty_{s/2}(0,t;b^{2-s}_{q,\infty}(\Om))}.
\end{array}
\end{equation}
Furthermore, using  $A_2^T A_2-A_1^T A_1= A_2^T \da A+ \da A^T A_1$,
\eq{4.45}, \eq{4.63n} and Lemma \ref{L4.7} (ii), we have
\begin{equation}
\label{E4.52}
\begin{array}{l}
\|\div [(A_2^T A_2-A_1^T A_1)\na v_1]\|_{L^\infty_{s/2}(0,t;b^{-s}_{q,\infty}(\Om))}
\ek
\qquad\leq \|\na (A_2^T \da A+ \da A^T A_1)\otimes\na v_1, (A_2^T
\da A+ \da A^T A_1)\otimes\na^2
v_1\|_{L^\infty_{s/2}(0,t;b^{-s}_{q,\infty}(\Om))} \ek
 \qquad\leq
c\big(\|\da A\|_{L^\infty(0,t; b^{1-s}_{q,\infty}(\Om))} \|\na
A_1\otimes\na v_1, \na A_2\otimes\na
v_1\|_{L^\infty_{s/2}(0,t;b^{-s}_{q,\infty}(\Om))} \ek\hspace{1cm}
+\|A_1, A_2\|_{L^\infty_{s/2}(0,t;b^{1-s}_{q,\infty}(\Om))} \|\na\da
A\otimes \na v_1,
  \da A\otimes \na^2 v_1\|_{L^\infty_{s/2}(0,t;b^{-s}_{q,\infty}(\Om))}\big)\ek
\qquad\leq c\big((\|\da A\|_{L^\infty(0,t; b^{1-s}_{q,\infty}(\Om))}
(\|\na A_1\otimes\na v_1, \na A_2\otimes\na v_1, \na^2
v_1\|_{L^\infty_{s/2}(0,t;b^{-s}_{q,\infty}(\Om))}\ek\hspace{1cm}
+\|A_1, A_2\|_{L^\infty_{s/2}(0,t;b^{1-s}_{q,\infty}(\Om))} \|\na^2
v_1\|_{L^\infty_{s/2}(0,t;b^{-s}_{q,\infty}(\Om))})\big)\ek
\qquad\leq c\|\da A\|_{L^\infty(0,t; b^{1-s}_{q,\infty}(\Om))}
\|A_1, A_2\|_{L^\infty_{s/2}(0,t;b^{1-s}_{q,\infty}(\Om))} \|\na^2
v_1\|_{L^\infty_{s/2}(0,t;b^{-s}_{q,\infty}(\Om))}\ek \qquad\leq
ct^{s/2}\|\da v\|_{L^\infty_{s/2}(0,t; b^{2-s}_{q,\infty}(\Om))}
\|v_1, v_2\|^2_{L^\infty_{s/2}(0,t;b^{2-s}_{q,\infty}(\Om))}.
\end{array}
\end{equation}

Thus, from \eq{4.57}-\eq{4.52} we have
$$\begin{array}{l}
\|\da f_2\|_{L^\infty_{s/2}(0,t;b^{-s}_{q,\infty}(\Om))}\ek
\qquad\leq ct^{s/2} \| v_1, v_2\|_{L^\infty_{s/2}(0,t;
b^{2-s}_{q,\infty}(\Om))} (\|v_1,v_2\|_{L^\infty_{s/2}(0,t;
b^{2-s}_{q,\infty}(\Om))}+1) \|\da v\|_{L^\infty_{s/2}(0,t;
b^{2-s}_{q,\infty}(\Om))},
\end{array}
$$
which together with \eq{4.50} yields
\begin{equation}
\label{E4.67} \|\da
F\|_{L^\infty_{s/2}(0,t;b^{-s}_{q,\infty}(\Om))}\ek \leq \eta_2(t)
\|\na^2\da v, \na\da Q\|_{L^\infty_{s/2}(0,t;
b^{-s}_{q,\infty}(\Om))}
\end{equation}
with some $\eta_2(t)$ such that $\eta_2(t)\ra 0$ as $t\ra +0$.

Finally, let us get estimate of $\|\da R_t\|$.
 Recall that $\da R=(\Id-A_2)\da v-\da A v_1$ and
hence $(\da R)_t=-A_{2t}\da v+(\Id-A_2)\da v_t-\da A_t v_1-\da A v_{1t}$.
Since we have
$$A_{it}(t,y)= -\Id+\sum_{k\geq 2}(-1)^kk C_i^{k-1}(t), i=1,2,$$
due to \eq{4.21}, it follows that
\begin{equation}
\label{E4.64}
\begin{array}{rcl}
\|A_{2t}\da v\|_{L^\infty_{s/2}(0,t; b^{-s}_{q,\infty}(\Om))} &\leq
&(1+\sum_{k\geq 2}kk_0^{k-1}\|\bar{C}_i(t)\|^{k-1}_{L^\infty(0,t)})
\|\da v\|_{L^\infty_{s/2}(0,t; b^{-s}_{q,\infty}(\Om))}\ek &\leq &c
t^{1-s/2}\|\da v\|_{L^\infty(0,t; b^{-s}_{q,\infty}(\Om))}
 \end{array}
 \end{equation}
 by repeatedly applying Lemma \ref{L4.7} (ii) in view of \eq{4.59}.
On the other hand, using \eq{4.54}, we have
\begin{equation}
\label{E4.65}
\begin{array}{rcl}
\|(\Id-A_2)\da v_t\|_{L^\infty_{s/2}(0,t; b^{-s}_{q,\infty}(\Om))}
&\leq & c \|\Id-A_2\|_{L^\infty(0,t; b^{1-s}_{q,\infty}(\Om))} \|\da
v_t\|_{L^\infty_{s/2}(0,t; b^{-s}_{q,\infty}(\Om))}\ek &\leq & c
\|\na v_2\|_{L^1(0,t; b^{1-s}_{q,\infty}(\Om))} \|\da
v_t\|_{L^\infty_{s/2}(0,t; b^{-s}_{q,\infty}(\Om))}\ek &\leq & c
t^{s/2}\|v_2\|_{L^\infty(0,t; b^{2-s}_{q,\infty}(\Om))} \|\da
v_t\|_{L^\infty_{s/2}(0,t; b^{-s}_{q,\infty}(\Om))}.
\end{array}
\end{equation}
Since $\da A_t= h_1 D\da v + h_3\int_0^t D\da v (\tau)\,d\tau, $
where $h_1$ is given by \eq{4.62} and
$$\begin{array}{l}
h_3:=h_{1t}=
  \sum_{k\geq 1}\sum_{j=0}^{k-1}\big(jC_1^{j-1}C_2^{k-1-j}D v_1+
  (k-1-j)C_1^{j}C_2^{k-2-j}D v_2\big),
\end{array}$$
we get, by applying Lemma \ref{L4.7} (iii) in view of \eq{4.59},
that
$$\|h_3\|_{L^\infty_{s/2}(0,t;b^{1-s}_{q,\infty}(\Om))}\leq
c\|v_1, v_2\|_{L^\infty_{s/2}(0,t; b^{2-s}_{q,\infty}(\Om))}.$$
Therefore, by Lemma \ref{L4.7} (ii) and \eq{4.60} we have
\begin{equation}
\label{E4.66}
\begin{array}{l}
\|\da A_t v_1\|_{L^\infty_{s/2}(0,t; b^{-s}_{q,\infty}(\Om))} \leq
c\|h_1\|_{L^\infty(0,t; {\cal M}(b^{-s}_{q,\infty}(\Om)))} \|D\da
v\otimes v_1\|_{L^\infty_{s/2}(0,t; b^{-s}_{q,\infty}(\Om))} \ek
 \hspace{4cm} +c\|h_3\|_{L^\infty_{s/2}(0,t;b^{1-s}_{q,\infty}(\Om))}\|D\da
v\|_{L^1(0,t; b^{-s}_{q,\infty}(\Om))}
 \ek
\leq ct^{s/2}\|v_1, v_2\|_{L^\infty_{s/2}(0,t;
b^{2-s}_{q,\infty}(\Om))} \big( \|D\da v\otimes
v_1\|_{L^\infty_{s/2}(0,t; b^{-s}_{q,\infty}(\Om))}+\|D\da
v\|_{L^\infty_{s/2}(0,t; b^{-s}_{q,\infty}(\Om))}\big)
 \ek
  \leq ct^{s/2}\|v_1, v_2\|_{L^\infty_{s/2}(0,t;
b^{2-s}_{q,\infty}(\Om))}\|\da v\|_{L^\infty_{s/2}(0,t;
b^{2-s}_{q,\infty}(\Om))}\big( \|v_1\|_{L^\infty(0,t;
b^{-s}_{q,\infty}(\Om))}+1\big).
\end{array}
\end{equation}
Finally, in view of the expression $\da A(t)=h_1(t)\int_0^t D\da
v\,d\tau$ and \eq{4.61}, we get by Lemma \ref{L4.7} (ii) that
\begin{equation}
\label{E4.68}
\begin{array}{rl}
\|\da A v_{1t}\|_{L^\infty_{s/2}(0,t; b^{-s}_{q,\infty}(\Om))} &\leq
c\|\da A\|_{L^\infty(0,t; b^{1-s}_{q,\infty}(\Om))}
\|v_{1t}\|_{L^\infty_{s/2}(0,t; b^{-s}_{q,\infty}(\Om))}\ek &\leq
ct^{s/2}\|\da v\|_{L^\infty_{s/2}(0,t; b^{2-s}_{q,\infty}(\Om))}
\|v_{1t}\|_{L^\infty_{s/2}(0,t; b^{-s}_{q,\infty}(\Om))}.
\end{array}
\end{equation}
Thus, from  \eq{4.64}-\eq{4.68} it follows that
\begin{equation}
\label{E4.69}
\|\da R_t\|_{L^\infty_{s/2}(0,t; b^{-s}_{q,\infty}(\Om))}
\leq \eta_3(t)( \|\da v\|_{L^\infty(0,t; b^0_{q,\infty}(\Om))}
+\|\da v\|_{L^\infty_{s/2}(0,t; b^{2-s}_{q,\infty}(\Om))}
\end{equation}
with some $\eta_3(t)$ converging to $0$ as $t\ra +0$.

Summarizing, we can conclude from \eq{4.36}, \eq{4.37}, \eq{4.49},
\eq{4.67} and \eq{4.69} that $\da v(t)=0, \da Q(t)=0$ for all $t\in
(0,T_1)$ with some $T_1>0$. Then, by standard continuation argument,
it  can be shown that $\da v(t)=0, \da Q(t)=0$ for all $t\in (0,T)$.

Now, the proof of the uniqueness part of Theorem \ref{T1.2} comes to
end. \qed

\newpage


%

\bigbreak \noindent {\bf Acknowledgments.}
P. Zhang is partially supported
by NSF of China under Grants   11371347 and 11688101, Morningside Center of Mathematics of The Chinese Academy of Sciences and innovation grant from National Center for
Mathematics and Interdisciplinary Sciences.


\begin{thebibliography}{99}


\bibitem{Abi07} H. Abidi, \'Equation de Navier-Stokes avec densit\'e et viscosit\'e variables
dans l'espace critique, Rev. Mat. Iberoam. 23 (2007), 537-586

\bibitem{AGZ12} H. Abidi, G. Gui and P. Zhang,
            On the wellposedness of three-dimensional inhomogeneous Navier-Stokes equations
in the critical spaces, Arch. Rational Mech. Anal. 204 (2012), 189-230


\bibitem{AGZ13} H. Abidi, G. Gui and P. Zhang,
               Well-posedness of 3-D inhomogeneous Navier-Stokes equations
               with highly oscillatory initial velocity field,
               J. Math. Pures Appl. 100 (2013) 166-203

\bibitem{Am00} H. Amann, On the strong solvability of the Navier-Stokes equations,
J. Math. Fluid Mech. 2 (2000), 16-98

\bibitem{Am95} H. Amann, Linear and Quasilinear Parabolic Problems, I.
Abstract Linear Theory, Birkh\"auser, 1995

\bibitem{BL77} J. Bergh and J. L\"ofstr\"om, Interpolation Spaces, Springer, 1977

\bibitem{Bo79} M. E. Bogovskii, Solution of the first boundary
            value problem for an equation of continuity of an
            incompressible medium. Dokl. Akad. Nauk SSSR. 248 (1979), 1037-1040.

\bibitem{BoSo90} W. Borchers and H. Sohr, On the equations $\text{rot} v = g$ and $\div u = f$
   with zero boundary conditions, Hokkaido Math. J. 19 (1990), 67-87.

\bibitem{Ca00} A. M. Caetano, Approximation by functions of compact support
in Besov-Triebel-Lizorkin spaces on irregular domains, Studia Math. 142 (2000), 47-63

\bibitem{DiLi89} R. J. DiPerna and P. L. Lions, Ordinary differential equations,
transport theory and Sobolev spaces, Invent. Math. 98 (1989), 511-547

\bibitem{Da03}  R. Danchin, Density-dependent incompressible viscous fluids
           in critical spaces, Proceedings of the Royal Society of Edinburgh,
                         133A, 1311-1334, 2003

\bibitem{Da06} R. Danchin, Density-Dependent incompressible fluids
            in bounded domains, J. Math. Fluid Mech. 8 (2006), 333-381

\bibitem{DaMu09} R. Danchin and P. B. Mucha, A critical functional framework
        for the inhomogeneous Navier-Stokes equations in the half-space,
         J. Func. Anal. 256 (2009), 881-927

\bibitem{DaMu12} R. Danchin and P. B. Mucha, A Lagrangian approach for the incompressible
Navier-stokes equations with variable density, Comm. Pure. Appl. Math.,
65 (2012), 1458-1480

\bibitem{DaMu13} R. Danchin and P. B. Mucha, Incompressible flows with piecewise constant
density, Arch. Rational Mech. Anal. 207 (2013), 991-1023

\bibitem{DaZh13} R. Danchin and P. Zhang, Inhomogeneous Navier-Stokes
equations in the half-space, with only bounded density, {\it J. Funct. Anal.} {\bf  267}  (2014),   2371-2436.

\bibitem{EvaGa92} L. C. Evans and R. F. Gariepy, Measure theory and fine properties of functions,
          CRC Press Inc, 1992

\bibitem{FaSo94} R. Farwig and H. Sohr, Generalized resolvent estimates for the
Stokes system in bounded and unbounded domains, J. Math. Soc. Japan
46 (1994), 607-643

\bibitem{Ga11} G. P. Galdi, An Introduction to the Mathematical Theory of the
Navier-Stokes Equations, Steady-State Problems, Second Edition, Springer, 2011

\bibitem{GHH05} M. Geissert, H. Heck and M. Hieber,
On the equation $\div u = g$ and Bogovskii's operator in Sobolev
spaces of negative order, Partial differential equations and
functional analysis, Oper. Theory Adv. Appl. 168 (2006), 113-121

\bibitem{GHP11} G. Gui, J. Huang and P. Zhang, Large global solutions to 3-D inhomogeneous
Navier-Stokes equations slowly varying in one variable, J. Func.
Anal. 261 (2011), 3181-3210

\bibitem{HPZ13} J. Huang, M. Paicu and P. Zhang, Global well-posedness
of incompressible inhomogeneous fluid systems with bounded density
or non-Lipschitz velocity, Arch. Rational Mech. Anal. 209 (2013),
631-682

\bibitem{Kaz74} A. V. Kazhikhov, Solvability of the initial-boundary value problem
for the equations of an inhomogeneous fluid. (Russian) Dokl. Akad. Nauk. SSSR 216 (1974), 1008-1010

\bibitem{Kat84} T. Kato, Strong $L^p$-solutions of the Navier-Stokes equations in
$\R^m$, with applications to weak solutions, Math. Z. 187 (1984), 471-480

\bibitem{LaSo75}  O. A. Ladyzhenskaya and V. A. Solonnikov, The unique solvability of
an initial-boundary value problem for viscous incompresible fluids. (Russian)
Boundary value problems of mathematical physics and related questions of the theory
of functions, 8 (1975), Zap. Nau\v{c}n. Sem. Leningrad. Otdel. Mat. Inst. Steklov (LOMI),
52, 52-109

\bibitem{PaZh12} M. Paicu and P. Zhang, Global solutions to the 3-D incompressible
           inhomogeneous Navier-Stokes system, J. Func. Anal. 262 (2012), 3556-3584

\bibitem{RiZh14} M.-H. Ri, P. Zhang and Z. Zhang, Global well-posedness for Navier-Stokes equations
with initial data in $B^0_{n,\infty}(\Om)$, J. Math. Fluid Mech.
18(2016), 102-131

\bibitem{Te77} R. Temam, Navier-Stokes equations, Noth-Holland, 1977

\bibitem{Tri83} H. Triebel, Interpolation Theory, Function Spaces, Differential Operators,
North Holland, 1983

\bibitem{Tri02} H. Triebel, Function spaces in Lipschitz domains
and on Lipschitz manifolds, Characteristic functions as a pointwise multiplers,
Rev. Mat. Complut. 15 (2002), 475-524

\bibitem{Ya10} A. Yagi, Abstract Parabolic Evolution Equations and their Applications,
        Springer, 2010


\end{thebibliography}
\end{document}